\documentclass[12pt]{amsart}
\usepackage{array}
\usepackage{longtable}
\usepackage{url}
\usepackage{enumerate}
\usepackage{amsmath}
\usepackage{amsthm}
\usepackage{amsfonts}
\usepackage{mathrsfs}
  \usepackage{ifthen}
 \newcommand{\mymarginpar}[1]{%
    \marginpar{\ifthenelse{\isodd{\arabic{page}}}{\flushleft #1}{\flushright #1}}}

 \renewcommand{\phi}{\varphi}

\newcommand{\IC}{\Bbb{C}}

 \newcommand{\IN}{\Bbb{N}}

 \newcommand{\IR}{\Bbb{R}}

 \newcommand{\IZ}{\Bbb{Z}}         
\newcommand{\FC}{\frak{C}}                  
 \newcommand{\eps}{{\varepsilon}}           
\newcommand{\CA}{\mathcal{A}}

\newcommand{\CC}{\mathcal{C}}
\newcommand{\CD}{\mathcal{D}}
\newcommand{\CE}{\mathcal{E}}
\newcommand{\CF}{\mathcal{F}}
\newcommand{\CS}{\mathcal{S}}
\newcommand{\CM}{\mathcal{M}}

\newcommand{\CQ}{\mathcal{Q}}

 \theoremstyle{plain} 
 \newtheorem{Theorem}{Theorem}[section]
 \newtheorem{Lemma}[Theorem]{Lemma}
 \newtheorem{Proposition}[Theorem]{Proposition}

 \theoremstyle{definition} 
 \newtheorem{Definition}[Theorem]{Definition}
 \newtheorem{Remark}[Theorem]{Remark}

\begin{document}

 \title{  Spectral triples for AF C*-algebras and metrics on the
   Cantor set. }
 
 \author{Cristina Antonescu and Erik Christensen}
 \address{Department of Mathematics \\University of California
   Riverside \\ 
Riverside, California 925074633 \\ USA. \\  \\Institut for Matematiske Fag University of Copenhagen\\
Universitetsparken 5\\
DK 2100 Copenhagen \O \\
Denmark }
 \email{cris@math.ucr.edu,  echris@math.ku.dk}
 \date{\today}


 \begin{abstract}
   
   An AF C*-algebra has a natural filtration as an increasing sequence
   of finite dimensional C*-algebras. We show that it is possible to
   construct a Dirac operator which relates to this filtration in a
   natural way and which will induce a metric for the weak*-topology
   on the state space of the algebra. In the particular case of a UHF
   C*-algebra the construction can be made in a way, which relates
   directly to the dimensions of the increasing sequence of
   subalgebras. It turns out that for AF C*-algebras once one has
   obtained a spectral triple, then the eigenvalues of that Dirac operator
   can be increased 
   arbitrarily without damaging the defining properties for a spectral
   triple. We have a obtained a version of an inverse to this result,
   by showing  that - under certain conditions, which are always true
   for the AF algebras we consider - 
   such a phenomenon  can only occur  for AF C*-algebras.

\smallskip \noindent
   The algebra of continuous functions on the Cantor set
   is an approximately finite dimensional C*-algebra and our
   investigations show, when applied to this algebra, that the
   proposed Dirac operators have good classical interpretations and
   lead to an, apparently, new way of constructing a representative
   for a Cantor set of any given Hausdorff dimension.  At the end of
   the paper we study the finite dimensional full matrix algebras over
   the complex numbers, $\CM_n$, and show that the operation of
   transposition on matrices yields a spectral triple which has the
   property that it's metric on the state space is exactly the
   norm distance. This result is then generalized to arbitrary
   unital C*-algebras.

 \end{abstract}

\maketitle

 \section{Introduction}
\noindent 
 Alain Connes has extended the notion of a compact metric space to the
 non commutative setting of C*-algebras
 and unbounded operators on Hilbert spaces, \cite{ Co2, Co3}. For a 
 compact, spin, Riemannian manifold $\mathcal{M}$, Connes has shown
 that the geodesic distance can be expressed in terms of an unbounded
 Fredholm module over the C*-algebra ${\rm C}(\mathcal{M})$, such that
 the
 distance between two points $p,q$ in $\mathcal{M}$ is obtained via
 the Dirac operator $D$ by the formula
 \begin{displaymath} d(p,q)= \sup \{\vert{a(p)-a(q)} \;\vert \; : \; a \in
   {\rm C}(\mathcal{M}), \Vert{[D,a]}\Vert\leq {1}\}.
\end{displaymath}

\noindent For a commutative unital C*-algebra ${\rm
  C}(X)$, the compact space $X$ 
 embeds naturally into the set of regular Borel
probability measures on the space $X$.
By Riesz' representation theorem the  latter space is the  weak*-compact
subset of the dual of ${\rm
  C}(X)$,  named the state space of the
algebra. For the case $X = \CM$, the notion
$d(p,q)$ 
can then be extended to
states on ${\rm C}(\mathcal{M})$ right away, using the same formula,
and in this way it is possible to obtain  a metric
for the weak*-topology on the state space of $ {\rm C}(\mathcal{M})$. 
 This construction does not rely on the commutativity of the
algebra ${\rm C}(\mathcal{M})$ but on the C*-algebra structure of  $
{\rm C}(\mathcal{M})$ and the existence of an operator like $D$, named
a {\em Dirac operator}. Given a C*-algebra $\CA$, the natural question is
then which properties an operator $D$ should have in order to deserve
such a name ? 
According to Connes, a candidate for a {\em Dirac
operator } must as a minimum have properties which he has defined in
the terms {\em unbounded
Fredholm module} and {\em spectral triple}.
\begin{Definition} \label{Fred}
Let $\mathcal{A}$ be a unital C*-algebra. An \textit{unbounded Fredholm module} 
$(H,D)$ over $\mathcal{A}$ is:
\begin{itemize}
\item[(i)]: a Hilbert space $H$ which is a left $\mathcal{A}$-module, that is, 
a Hilbert space $H$ and a *-representation of $\mathcal{A}$ on $H$;

\smallskip
\item[(ii)]: an unbounded, self-adjoint operator $D$ on $H$ such that the set
\begin{displaymath}
 \{a\in\mathcal{A}:[D,a] \text{ is densely defined and extends to a
  bounded operator on } H \}
\end{displaymath}
is norm dense in $\mathcal{A}$;

\smallskip
\item[(iii)]: $(I+D^{2})^{-1}$ is a compact operator.
\end{itemize}

\smallskip
\noindent
The triple $(\mathcal{A},H,D)$
with the above description is called \textit{a spectral triple}.
\end{Definition}

\smallskip
\noindent
Condition (iii) is
quite often strengthened in the way  that $D$ is said to be finitely
summable or $p$-\textit{summable}, 
\cite{Co2, Co3}, if for some $p > 0$
\begin{displaymath}
\text{trace}\left ( \left ( I +D^{2} \right )^{-p/2} \right )<\infty
\end{displaymath}

\medskip

\noindent
Given a spectral triple  $(\mathcal{A},H,D)$,  
one can then introduce a pseudo-metric on the state 
space $\CS(\CA)$ of $\mathcal{A}$ by the formula
\begin{displaymath}
\forall \phi, \psi \in \CS(\CA):\; d(\phi,\psi)= \sup\{\vert\phi(a)-\psi(a)\vert: a\in \mathcal{A},\Vert[D,a]\Vert\leq1\}.
\end{displaymath}

\noindent
We use the term \textit{pseudo-metric} because it is not clear that
$d(\varphi,\psi)<\infty$ for all pairs, but  the other axioms 
of a metric are fulfilled.

\noindent
Marc A. Rieffel has studied several aspects of this extension of the concept
of a compact metric space to the framework of C*-algebras, and he has
obtained a lot of results \cite{Ri2, Ri3, Ri4, Ri5, Ri6}. Among the questions he has
dealt with, we have been most attracted by the 
 one which asks whether a spectral triple
will induce a metric for the weak*-topology on the state space.
If the metric topology coincides with the weak*-topology on the state space,
then the metric topology should give the state space a finite diameter, since 
the state space is compact for the weak*-topology.
A nice  characterization of when the metric is bounded  
on the state space and furthermore when it induces the weak*-topology
on this space
was given by Rieffel, \cite{Ri2}, and Pavlovi\'c, \cite{Pav}. 
This characterization reads:
\begin{Theorem} \label{cmpmet}
Let $(H,D)$ be an unbounded Fredholm module over a unital C*-algebra
$\mathcal{A}$, and let the pseudo-metric $d$ on $\mathcal{S}(\mathcal{A})$ be 
defined by the formula:
\begin{displaymath}
d(\phi,\psi)= \sup\{\vert\phi(a)-\psi(a)\vert: a\in \mathcal{A},\Vert[D,a]\Vert\leq1\}.
\end{displaymath}
for $\phi,\psi\in \mathcal{S}(\mathcal{A})$ then
\begin{itemize}
\item[(i)] $d$ is a bounded metric on $\mathcal{S}(\mathcal{A})$ 
if and only if 
\begin{displaymath}
\{a\in\mathcal{A}:\Vert [D,a]\Vert\leq 1\}
\end{displaymath}
has a bounded image in the quotient space $\mathcal{A}/\Bbb{C}1$, equipped
with the quotient norm.\newline
\item[(ii)] the metric topology coincides with the weak*-topology if and only if the set 
\begin{displaymath}
\{a\in\mathcal{A}:\Vert [D,a]\Vert\leq 1\}
\end{displaymath}
has a precompact image in the quotient space $\mathcal{A}/\Bbb{C}1$,
equipped with the quotient norm.  
\end{itemize}
\end{Theorem}

\noindent
We will recall some classical results, which serve as motivation for
the definitions given in the non commutative case.
 The compact manifold $\Bbb{T}$, i.e. the unit circle, has a
 canonical differential operator and is also a compact group. This
 group is  the dual
 of the integers and the integers has a canonical length function, which,
 via the Fourier transform, 
 has close relations to the differentiation operator on the circle. 
While respecting all these classical structures,  
 Connes  considered in \cite{Co2} 
a discrete group $G$ endowed with a \textit{length function}  
$\ell :G\to\Bbb{R}_{+}$
In analogy with the situation for $C(\Bbb{T})$, or the
$2\pi$-periodic
 functions on $\IR$ where $\frac{1}{i}\frac{d}{dt}e^{int} = ne^{int}$,
 Connes then 
defined a Dirac operator $D$ on $\ell^{2}(G)$ by
$(D\xi)(g)=\ell(g)\xi(g)$ and he
 proved, that if the length function $\ell$ is a \textit{proper length function}, i.e. 
$\ell^{-1}([0,c])$ is finite for each $c\in \Bbb{R}_{+}$,
then $(\ell^{2}(G),D)$ is an unbounded Fredholm module for
$C^{*}_{r}(G)$. It is a remarkable result of \cite{Co2}, that 
 only amenable
discrete groups can have a  $p-$summable Fredholm module. This result
is extended by Voiculescu in \cite{Vo}, sections 4 and 5. The common features
of these results are that certain boundedness properties of  a 
spectral triple $(\CA, H, D) $  imply that $\CA$ has a faithful 
 trace state
which extends to an $\CA$ invariant state on $B(H)$.  This in turn
implies \cite{Co1} that the representation of $\CA$ on $H$ must be
hyperfinite, and in case of a reduced group C*-algebra of a discrete
group, it therefore
follows that the group  must be amenable. In particular these results
tell that for certain spectral triples $(\CA, H, D)$ involving a non nuclear
C$^*$-algebras it is not possible to perturb $D$ very much without
destroying the properties of a spectral triple. It is one of the main
insights in the present investigation that for some of the  spectral triples
$(\CA, H, D )$ one can construct for  an approximately finite
dimensional C$^*$-algebras there is an abundance of possible
perturbations of $D$ which will still give a spectral triple. 

\medskip
\noindent In \cite{Co2, Co3, Vo, OR} and many more places the concept {\em 
filtration of a C$^*$-algebra }  plays an important role in the
investigation of spectral triples. The reason is of course that the
filtration quite often induces a natural candidate for a spectral
triple. The present authors were inspired by this and then wanted to
see what this line of investigation can yield for approximately finite
dimensional, or AF, C$^*$-algebras.

\smallskip
\noindent
 An
AF C*-algebra $\mathcal{A}$ has a natural filtration since
$\mathcal{A}$, by definition, is the norm closure of an increasing
sequence of finite dimensional C*-algebras
$(\mathcal{A}_n)_{n \in \Bbb{N}_{0}}$. These algebras were
studied first by Bratteli, \cite{Br}. In this paper we will only 
consider unital  AF C*-algebras which have  a faithful 
state, and we will therefore always assume that
$\mathcal{A}_0=\Bbb{C}I_{\mathcal{A}}$. 
Based on the GNS-representation coming from  a faithful state
and the given increasing sequence of subalgebras, we show that it is
possible to construct an unbounded Fredholm module over this C*-algebra
in much the same way as it was done by Connes for discrete
groups.  To verify the agreement between the induced metric topology
and the weak*-topology on the state space we follow ideas of the
same type as was used in the group C*-algebra case, 
\cite{AC, Co2, Co3, OR, Ri4}. This means that
we try to obtain estimates of the norm of an element $a$ by careful 
estimation of the norm of some parts  of the
commutator $[D,a]$. The hole point in these computations is to show
that there exists a $D$ such that the set $\CD = \{a \in \CA \, : \,
\|[D,a]\| \leq 1 \, \}$ is so big that it separates the states of
$\CA$ and so small that it has precompact image in $\CA/\IC I$.
It is not
difficult to see that such a Dirac operator must exist, and further
that there is a lot of freedom in the choice of the eigenvalues for
such a $D$. Especially for the Dirac operators obtained,  it turns out
that a certain minimal growth in the eigenvalues is needed in order to
get the right properties, but once this level is attained, it is
possible to increase the numerical values of eigenvalues arbitrarily,
without destroying  the topological properties of the spectral
triple. This implies that for any $p > 0$ it is
possible to construct a $p-$summable Fredholm module, so it is not
possible to assigne a dimension other than 0 to a unital AF C*-algebra via
spectral triples.

\smallskip
\noindent
 As mentioned above this is very much in contrast to results by Connes 
and Voiculescu  \cite{Co2, Vo}. On the other hand it
 was suggested to us by Connes that the possibility to
 increase the numerical values of the  eigenvalues of a Dirac
operator can only happen for AF C*-algebras. We have included a
theorem which confirms that conjecture if the eigenspaces of the Dirac
operator is fixed.

\smallskip
\noindent
In the special case of a UHF
C*-algebra, the Bratteli diagram is easy to analyze and the unique trace
state is faithful, so it is quite easy to 
give a {\em natural } description of a spectral triple with the right
properties. 

\smallskip
\noindent
Another special type of AF C*-algebras are the commutative ones. The
approximately finite dimensionality implies here, that such an algebra
consists of  the
continuous complex functions on a totally disconnected compact space.  
It is well known that the algebra of continuous functions on the 
standard Cantor subset of the unit interval is an approximately finite
dimensional C*-algebra, so we have tried to see what the spectral
triples could look like  in this case. We are not the first ones who
try to apply the non commutative tools on this commutative
algebra. We have had the opportunity to see some notes by Connes
\cite{Co4}, where he constructs a spectral triple for the algebra of
continuous functions on the Cantor set. His construction is different
from ours and  very accurate in its reflection of the geometrical
properties of the usual Cantor set which is obtained by successive
cuttings of open intervals from the unit interval. In particular
Connes spectral triple  
makes it possible to recover the metric inherited from $\IR$ exactly
and he can find the differentiation operator too. Our emphasis is to
see how a  general construction which works for any AF C*-algebra will
work in this special case. This has the effect that the module
we propose is quite different from the one Connes has constructed. 
We  try, later, 
shortly to describe the major difference between the 2 types of
modules. Our spectral triple will not give exact geometric data for
the middle third Cantor set, but it will in this case 
induce a metric equivalent to the one coming from the standard
embedding of the Cantor set in $\IR$. On the other hand  
we can see that there are several other possible choices of spectral
triples and 
we have found a {\em ``natural''}
 family $(D_{\gamma})_{\gamma \in ]0,1[}$ of Dirac
operators, which all act on the Hilbert space coming from the
standard representation of the Cantor algebra. 
It turns out that for any $\gamma
$ in $]0,1[$ the corresponding Dirac operator will yield a metric on
the Cantor set, such that this compact  metric space, say
$\CC_{\gamma}$, 
 will have
Hausdorff dimension $\log 2 / (- \log \gamma)$. Having this, we
searched the literature for a unified representation theory for
Cantor sets of any positive dimension. We have not found a general
theory, but we found some examples \cite{ FY}, \cite{Ke}.
These examples  have only a small   connection to our spectral
triples, so we have constructed a
  family of compact subsets of $\ell^1(\IN,
\IR) $ which serves our purpose. Each of theses spaces is bi-Lipschitz
equivalent to {\em a generalized Cantor set } in some space $\IR^e$.
Following a suggestion by Marc Rieffel we have computed an upper bound for
the Gromov--Hausdorff distance between two such compact metric spaces
$\CC_{\gamma} \text{ and } \CC_{\mu}$.
 Finally this investigation showed us how to construct 
 a compact metric space which for any $\gamma \in
]0,1[$ contains a subset which is bi-Lipschitz equivalent to
$\CC_{\gamma}$. 

\smallskip
\noindent
A full matrix algebra $\CM_n(\IC)$ is a special UHF C*-algebra and it
has a very special compact metric on the state space, namely the one 
 induced by the norm. We
show that for $\CM_n(\IC)$ acting on itself with respect to the
trace state and the Dirac operator given by
transposition on   $\CM_n(\IC)$, we get  a spectral triple such
that the norm distance on the state space is the metric induced by
the Fredholm module.
\noindent
Some more investigations into this construction show that it is possible
to extend the properties of the transposition operator to the setting
of a general C*-algebra $\CA$. Having this, it turns out that the metric
induced by the norm on $\CS(\CA)$ can always be obtained via a Dirac
operator which is a self adjoint unitary.

\section{AF C*-algebras} \label{AF}
\noindent
We consider now the case of AF C*-algebras.  Let $\mathcal{A}$ be a
unital AF C*-algebra, such that $\CA$ is the norm closure of an
increasing sequence $(\CA_n)_{n \in \IN_0}$ of finite dimensional
C*-algebras.
 As mentioned above we
will always stick to the case of a unital C*-algebra $\CA$ and assume
that $\mathcal{A}_0=\Bbb{C}I_{\mathcal{A}}$. Further we assume that
$\CA$ has a faithful state, which we denote $\tau$.  For a
natural number $k$ we will let $\CM_k$ denote the algebra of complex
$k \times k$ matrices. Then for $n\geq 1$ each $\mathcal{A}_{n}$ is
isomorphic to a sum of full matrix algebras
$\mathcal{M}_{m_1}\oplus\ldots\oplus\mathcal{M}_{m_{k}}$ and it embeds
into $\mathcal{A}_{n+1}$ such that the unit of $\CA$ always is the
unit of each algebra in the sequence.  The GNS Hilbert
space $H$ of $\mathcal{A}$ is just the completion of pre Hilbert space
$\mathcal{A}$ equipped with the inner product $(a,b)=\tau(b^{*}a)$,
and the GNS representation $\pi$ is just given by left multiplication.
In order to avoid the writing of too many $\pi's$ we will assume that
$\CA$ is a subalgebra of $B(H)$ and that $\xi $ is a unit vector in $H$
which is separating and cyclic for $\CA$ and further has the property
that the vector state $\omega{_\xi} $ equals $\tau$.  Since $\xi$ is
separating, the mapping $\eta : \CA \to H $ given by $\eta(a) = a\xi$
induces - for each $n$ in $\IN$ - a bijective linear homeomorphism of
the algebra $\CA_n $ onto a finite dimensional subspace - say $H_n$ of
$H$.  The corresponding growing sequence of orthogonal projections
from $H$ onto $H_n$ is denoted $(P_n)$.  The basic idea in the 
construction of  a Dirac operator is, to let it have its eigenspaces
equal  to the sequence of differences $H_{n} \ominus H_{n-1}$. We do
therefore define a sequence $(F_n)_{n\in \IN}$ of pairwise orthogonal finite
dimensional subspaces of $H$ by $F_0 = H_0$ and $F_n = H_{n} \ominus
H_{n-1}$ for $n\geq 1$. The corresponding sequence of pairwise
orthogonal projections is denoted $(Q_n)_{n\in \IN}$, so $Q_0 = P_0$ and $Q_n =
P_n - P_{n-1}$ for $n\geq 1$. We will, in the following theorem,
  show that it is
possible to find a sequence of positive reals $(\alpha_n)_{n \in
  \IN_0}$ such that the operator $D$ given by

\begin{equation*} \label{Diracdef}
D=\overset{\infty}{\underset{n=1}{\sum}}{\alpha_nQ_{n}}
\end{equation*}
can serve as Dirac operator for a reasonable spectral triple.
\noindent Remark that $\alpha_0$ is not needed in the
description above so we will therefore fix $\alpha_0 = 0$. In
particular this means that $D\xi = 0$. Based on all this notation we
can now formulate the main result of this section which in combination
with Theorem \ref{cmpmet} shows that the metric induced by $D$ on the
state space of $\CA$ will be a metric for the weak*-topology on the
state space.

\begin{Theorem} \label{DAF}
  Let $\CA$ be an infinite dimensional unital AF C*-algebra acting on
  a Hilbert space $H$ and let $\xi $ be a separating and cyclic, unit
 vector for $\CA$. 

\begin{itemize}
\item[(i)] There exists a sequence of real numbers
  $(\alpha_n)_{n \in \IN_0}$ such that $\alpha_0 = 0$ and  with the
  notation introduced above, the unbounded
  selfadjoint operator
  $D=\overset{\infty}{\underset{n=1}{\sum}}{\alpha_nQ_{n}}$ on $H$ has
  the property that the set
\begin{displaymath}
\CD = \{a = a^* \, \in \CA \, : \, \|\,[D,a]\, \| \leq 1  \}
\end{displaymath}
separates the states and 
has a precompact image in the quotient space $\mathcal{A}/\Bbb{C}1$,
equipped with the quotient norm.  Further the metric which $\CD$
induces on the state space generates the weak*-topology.
\item[(ii)] Given any $p > 0$ it is possible to choose the sequence $(\alpha_n)_{n \in
    \IN_0}$, such that the Fredholm module is $p$-summable.
\end{itemize}
\end{Theorem}   
\begin{proof}
We remind
the reader of the fact that the algebras $\mathcal{A}_n$ are finite dimensional
 and as such, they are all  
complemented subspaces of $\CA$. In the special case where  $\xi$ is a
trace vector it is well known how to construct a completely positive
projection of $\CA$ onto $\CA_n$, \cite{Ch}.
 In the general case with $\xi $ just
separating and cyclic, the method can be mimicked and  we can 
 define a continuous projection, say $\pi_n$
of $\CA$ onto $\CA_n$ by 

\begin{displaymath}
\forall \, a \,   \in \CA: \quad \pi_n(a)\,  :=  \,
\eta^{-1}(P_n\eta(a)).
\end{displaymath} 
 This definition
has one fundamental consequence upon which we shall build our
arguments, namely, that when we consider both the operator algebra
structure and the Hilbert space structure on $\CA$ simultaneously we
get
 \begin{displaymath}
\forall a \in \CA: \; \pi_n(a)\xi = P_na\xi.
\end{displaymath}

\noindent
In order to prove that the set $\CD$ has the properties stated we
start by studying the domain of definition for the unbounded
derivation 
$\delta(a  ) = [D,a]$. Hence we think that a sequence $(\alpha_n)_{n
  \in \IN_0} $ with $\alpha_0 = 0$ is given and we will show that 
for any $n \in \IN$ and any $a \in \CA_n$ the
commutator $[D, a]$ is densely defined and bounded. This shows that 
the union of multiples of $\CD$ given as 
$\cup_{n \in \IN}n\CD$ is dense in the self adjoint part of 
$\CA$ and consequently
must separate the states on $\CA$. So let us fix an $n$ in $\IN$,
a self adjoint   $ a $ in $ \CA_n$ and let $m$
in $\IN$ be chosen such that $m > n$. Since $a\CA_m \subset \CA_m$ and
 $a\CA_{m-1} \subset \CA_{m-1}$ we find that both of the projections
 $P_m$
 and $P_{m-1}$
commute with $a$ and consequently $a $ commutes with $Q_m = P_m -
P_{m-1}$. This means that for the closure of the commutator $[D,a]$
we get  
\begin{displaymath}
\forall n \in \IN \, \forall a \in \CA_n: \; \text{closure}( \, [D,a] \,)
= \sum_{i=1}^n \alpha_i [Q_i,a].
\end{displaymath}
\noindent
Since the sum above is finite, some positive multiple of 
$\, a \, $must be in
$\CD$. It is important to notice that this is true  for any sequence 
$(\alpha_n)$, and the   statement that $\CD$ separates the
states will be independent of the actual size of the eigenvalues
$\alpha_n$.

\noindent
The main ingredient  of our construction is that we show that there
exists a sequence of positive reals $(c_n)_{n \in  \IN}$ which only
depends on the sequence $(\CA_n)$ and the vector $\xi$ such that
for any
$ a \in \CD$ -  where now $\CD$ depends on the actual values of $\alpha_n$-
we have 

\begin{displaymath} 
\forall n \in \IN_0 \, \forall a \in \CD: \;  \Vert
\pi_{n+1}(a)-\pi_{n}(a)\Vert \; \leq \; \frac{c_{n+1}}{\alpha_{n+1}}
\end{displaymath}
\noindent
When this relation is established it follows easily that we can get
very nice convergence estimates by choosing the eigenvalues $\alpha_n$
sufficiently big.
\noindent
We will now describe how the sequence $(c_n) _{n \in \IN} $ is
obtained. For any $n \in \IN_0$  the seminorms $ a \to \Vert
\pi_{n+1}(a)-\pi_{n}(a)\Vert_{\CA}$ and  $ a \to  \| Q_{n+1}a\xi \|_H$ are
equivalent, since $Q_{n+1} $ is of finite dimension and 
$\xi$ is a separating vector. 
 Consequently there exists a positive real $c_{n+1} \geq 1$ such that
\begin{displaymath}
\forall n \in \IN_0 \, \forall a \in \CA: \;
\|\pi_{n+1}(a)-\pi_{n}(a)\|\; \leq \;   c_{n+1} \| Q_{n+1}a\xi \|.
\end{displaymath} 

\smallskip
\noindent
We will now again suppose that a sequence $(\alpha_n)$ is given such
that $\alpha_0 = 0$,
 fix  an arbitrary $n\in\Bbb{N}_0$ and an $a$ in $\CD$, then since
 $D\xi = 0$
we get a series of estimates
\begin{alignat*}{2} \label{ineqinAF}
1 \quad & \geq \quad \|[D,a]\|  \quad & \geq &  \quad \|Q_{n+1}[D,a]Q_0\| \\
  &  = \quad \|\alpha_{n+1}Q_{n+1}a \xi\|  \quad & \geq & \quad
  \frac{|\alpha_{n+1}|}{c_{n+1}} \|\pi_{n+1}(a) - \pi_n(a)\|.\\
\end{alignat*}
\noindent
In conclusion for any $n\in\Bbb{N}_0$ and any $a\in\mathcal{D}$ we have 
\begin{equation*} \label{pindif}
\Vert \pi_{n+1}(a)-\pi_{n}(a)\Vert\leq \frac{c_{n+1}}{|\alpha_{n+1}|}. 
\end{equation*}
We will now make the choice of the elements in the sequence
$(\alpha_n)_{n \in \IN_0}$ such that the sequence of fractions
$(\frac{c_{n+1}}{|\alpha_{n+1}|})_{n \in \IN_0}$ is summable. 
Let us then consider an absolutely  convergent series $
\sum_{\IN}\beta_n$  of positive reals and for  $n \in \IN $ we define 
$\alpha_n = \beta_n^{-1} c_n$  
\begin{align*} \label{aandif} 
\forall a \in \CD \forall n \in \IN \forall k \in \IN:\;  \|\pi_{n+k}(
a) - \pi_n(a) \| & \leq
\sum_{j= 1}^{k} 
\|\pi_{(n+j)}(a) - \pi_{n+j-1}(a)\| \\
& \leq \sum_{j=1}^{\infty} 
\frac{c_{n+j}}{\alpha_{n+j}} = \sum_{j=1}^{\infty}\beta_{n+j}.  
\end{align*}
From this inequality we first deduce that the sequence $(\pi_n(a))_{n
  \in \IN} $ is a Cauchy sequence in  $\CA$ and hence convergent. Let $b $
denote the limit point for this sequence, then 
\begin{displaymath}
b\xi = \lim \pi_n(a)\xi = \lim P_n a\xi = a\xi,
\end{displaymath}
and since $\xi$ is separating $b =a$ and we get from above the
following estimates;
\begin{equation*}  
\forall a \in \CD \forall n \in \IN :\;  \|a - \pi_n(a) \| \leq
\sum_{j=1}^{\infty}\beta_{n+j}.  
\end{equation*}

\noindent
For $n =0$ the inequality above gives
\begin{displaymath}
\forall a \in \CD: \; \| a - \omega_{\xi}(a)I \| \; = \; 
\| a - \pi_0(a) \| \; \leq \; \sum_{j=1}^{\infty} \beta_j < \infty. 
\end{displaymath}
In particular this shows that for any $a$ in $\CD$, the norm of $a +
 \IC  I$ in the quotient space is at most $  \sum \beta_k$, so
$\CD/\IC I$ is bounded.

\noindent The inequality also yields the precompactness of $\CD/\IC I$
right away. Let $\eps > 0 $ be given and  let $ n $ in $\IN$ be chosen
such that $\sum_{j \geq n}\beta_{j+1} < \frac{\eps}{2}$, then
\begin{displaymath}
\forall a \in \CD: \; \| a - \pi_n(a) \| \; \leq \; \frac{\eps}{2} \;
\text{ and } \; \| \pi_n(a) - \omega_{\xi}(a) I \| \; \leq \;
\sum \beta_j. 
\end{displaymath}
Since $ \CA_n$ is finite dimensional; a  closed ball in this space of radius 
 $\sum \beta_k$ is norm compact and can be covered by finite number
 of balls of radius $\frac{\eps}{2}$. This shows that the
 set
\begin{displaymath} 
\{ a - \omega_{\xi}(a) I \; : \; a \in \CD \}
\end{displaymath}
can be covered by
 a finite number of balls of radius $\eps$ and hence this set is precompact in
 $\CA$ and consequently $\CD / \IC I $ is precompact in $\CA/ \IC I$.

\medskip
\noindent
Let $p > 0 $ be given. 
With respect to the $p-$summability of the above Fredholm module we
may assume that the sequence of algebras satisfies $\CA_{n+1} \neq
\CA_n$. If this was not so, then $Q_{n+1} = 0$ and the repetition of
$\CA_n$ would not have any impact on the Dirac operator. We will
therefore assume that for all $n$ in $\IN_0$ we have
$\text{dim}\CA_{n+1}\, > \, \text{dim}\CA_n$. Since $\CA_0 $ is one
dimensional we get then the  rough estimates 
\begin{displaymath}
\forall n \in \IN_0: \quad  \text{dim}\CA_n \geq n+1 \text{ and  } \text{dim}\CA_{n+1}\, > \, \text{dim}\CA_n .
\end{displaymath}

\noindent
For the given $p$ we define $t = \max\{\, 2, \, \frac{3}{p}\,\}$, then
for $\beta_n = (\text{dim}\CA_n)^{-t}$ we have $\sum \beta_n < \infty$
since 
\begin{displaymath}
\sum_{n=1}^{\infty} \beta_n \;   \leq \; \sum_{n=1}^{\infty}
(\text{dim}\CA_n)^{-2} \; \leq \; \sum_{n=1}^{\infty} (n+1)^{-2} \;
\leq \; 1.
\end{displaymath}
\noindent
The $p-$summability follows in the same way since we know for all $n$
we have 
$c_n \geq 1 $. 
\begin{align*}
\text{tr}( ( I + D^2)^{-\frac{p}{2}})\, &= \, 1 + \sum_{n=1}^{\infty} ( 1 +
(\text{dim} \CA_n)^{2t}c_n^2)^{-\frac{p}{2}} (\text{dim}\CA_n -
\text{dim}\CA_{n-1} ) \\ 
&  \leq \, 1 +  \sum_{n=1}^{\infty}
(\text{dim}\CA_n)^{(1-pt)} \\
& \leq \, 1 +
\sum_{n=1}^{\infty}(n+1)^{-2} \, \leq \, 2.
\end{align*}
\end{proof}

\begin{Remark}
It may, in the first place, seem plausible that the construction given
just above should be
applicable in a wider setting than just the one of AF C*-algebras. We
have - of course - tried to follow such a path, and realized that -
at least for us -  it is not 
easy to get much further along this road.
Suppose for a moment that the elements in the filtration
$(\CA_n)$ 
are no longer algebras but just  finite dimensional
subspaces such that $\CA_n \CA_m \subset \CA_{m+n}$, then the spaces 
$H_m = \CA_m \xi$ can not be expected to be invariant for operators
$a$ in $\CA_n$ when $ n \leq m$.
In this case, the boundedness of the commutator $[D,a] $, for $a \in
\CA_n$
 can not be
established easily unless the sequence $(\alpha_n)_{n \in \IN_0}$
is the very special one given as $\alpha_n = n$. To realize that it is,
in general, not possible just to increase the eigenvalues arbitrarily
we refer, again, to Connes' \cite{Co2} and Voiculescu's results
\cite{Vo}. In particular they show  that the reduced group
C*-algebra of a  non amenable discrete group
can not have a finitely summable Fredholm module.  
On the other hand for AF C*-algebras there is no upper limit to the
growth of the eigenvalues of $|D|$. This is not only remarkable when
compared to the just mentioned results of Connes and Voiculescu,
 but also quite 
opposite to classical results for commutative C*-algebras associated
to compact smooth manifolds. In the commutative world an AF C*-algebra
is the continuous functions on a totally disconnected compact space
and the present general C*-algebraic result fits well with this fact.
Yet another aspect is discussed in \cite{Vo}, namely the possibility
of having a slow growth of the dimensions dim$\CA_n$ of the elements
in the filtration. We have not been able to obtain results for
spectral triples related to a filtration with slow growth, but it
seems likely that such an assumption might have non
commutative geometrical consequences. 

\smallskip
\noindent
In a presentation of Theorem
2.1 to an audience containing A. Connes, he suggested that this freedom
in the choice of unbounded Fredholm modules might be a characteristic
property for AF C*-algebras. We have tried to prove this conjecture,
but in vain for this sort of generality. On the other hand 
  we have found a more
restrictive property - which of course always holds for an AF
C*-algebra -  and  we can show that algebras which have this property
must be AF C*-algebras.
\noindent
Our difficulty in proving a general result lies in the problems
involved in comparing Dirac-operators connected to different
representations 
and/or with different spectral projections. If we fix the
representation, and   the spectral projections of $D$ and further impose  an
extra condition on some dense subset of the algebra,  then we can obtain a
result of the type conjectured by Connes.

\end{Remark}

\begin{Theorem} \label{AFD}
  Let $\CA$ be a C*-algebra acting on
  a Hilbert space $H$,  $\xi $ be a separating unit
  vector for $\CA$ and $(Q_n)_{n \in \IN_0}$ a sequence of pairwise
  orthogonal finite-dimensional projections with sum $I$ such that 
$ Q_0 \xi = \xi$.
   For any
  sequence of real numbers $(\lambda_n)_{n \in \IN_0}$ such that
  $|\lambda_n| \to \infty \text{ for } n \to \infty$  the symbol
  $D_{(\lambda_n)}$ shall denote the closed self adjoint operator which
  formally can be written as $\sum \lambda_n Q_n$. The common domain of
  definition, $ \rm{span}$  $( \cup Q_nH)$ for all the operators $D_{(\lambda_n)}$  is denoted $\CD_0$. 
  
\noindent If $\CA$ contains a dense subset $\CS$ such that for any $s$
in $\CS$ and any $D_{(\lambda_n)}$ the commutator $[D_{(\lambda_n)},
s]$ is  defined and bounded on  $\CD_0$ then $\CA$ is an AF
C*-algebra.
\end{Theorem}   
\begin{proof}

The sequence $(Q_n)_{n \in \IN_0}$ induces a matrix description of the
operators on $H$ so we will define $H_n = Q_nH$ and for an operator
$x$ in $B(H)$ we will write $x = (x_{ij})$ such that $x_{ij} $ is an
operator in $B(H_j,H_i)$ given by $x_{ij}= Q_ix|H_j$. For any bounded
operator $x$ and any sequence $(\lambda_n)$ we can then formally describe 
the commutator $[D_{(\lambda_n)}, x]$ by 
\begin{displaymath}
\forall i, j \in \IN_0:\;[D_{(\lambda_n)}, x]_{ij} = (\lambda_i -
\lambda_j)x_{ij}.
\end{displaymath}
\noindent
This may be nothing but a formal infinite matrix, but if we know that
the commutator is defined and bounded on $\CD_0$ then the  closure of
the commutator will be a bounded
operator whose matrix is the one just described.

\smallskip
\noindent
Given the assumptions on the elements in $\CS$ we may then work on
the matrix representations of their commutators
with various $D_{(\lambda_n)}$'s without worrying  about the domain of
definition for the commutators. Our first aim is to prove that for
elements in $\CS$ all, but finitely many,  of the matrix entries outside
the main diagonal vanish.
Let us then assume that there is an element $s$ in $\CS$ which has
infinitely many non zero matrix entries  outside the main diagonal,   
  and let
$(s_{{i_k}{j_k}})_{k \in \IN}$ be an infinite 
 sequence of non vanishing entries
such that for all $k$,  $i_k \neq j_k$.
Without loss of generality we may assume that the sequence
$(n_k)_{k\in \IN}$
of natural
numbers given by $ n_0 = 0 \text{ and } n_k := \max\{i_k, j_k\}$ for
$k \geq 1 $  is strictly increasing. We will then inductively over $ k
\in \IN$ define a
sequence $(\lambda_n)_{n \in \IN_0} $ - which will yield the contradiction - by

\begin{align*}
\lambda_0 &= 0 \\
\lambda_n &= \lambda_{n_{k-1}} \text{ if } n_{k-1} < n < n_k \\
\lambda_{n_k} &= \lambda_{n_{k-1}} + k + \frac{k}{\|s_{i_kj_k}\|}  
\end{align*}

\noindent
For this sequence $(\lambda_n)$ and $k$ in $\IN$ we want to estimate
the norm of the matrix entry  
$[D_{(\lambda_n)}, s]_{i_kj_k} = (\lambda_{i_k} -
\lambda_{j_k})s_{i_kj_k}$. Since $i_k \neq j_k$  the norm 
will satisfy $ \|[D_{(\lambda_n)},
s]_{i_kj_k} \| \, \geq  \, k$,  so $[D_{(\lambda_n)}, s]$ is not
bounded on $\CD_0$ and any element $s $ in $\CS$ has only 
finitely many
non zero entries outside the main diagonal 
in it's  matrix $(s_{ij})$. In particular this means
that for each $s$ in $\CS$ there exists a natural number $n_s$ such
that for any natural number $n > n_s$ 
the projection $Q_n$ will commute with  $s$.
In order to get our skeleton of finite dimensional C*-subalgebras of
$\CA$ we then define an increasing  sequence of finite dimensional 
projections $(P_n)_{n \in \IN_0}$ in $B(H)$ by 
\begin{displaymath}
\forall n \in \IN_0 \quad P_n \, := \, \sum_{i =0}^n Q_i.
\end{displaymath} 

\noindent
For each $n \in \IN_0 $ we can then define an increasing sequence of unital
C*-subalgebras $\CA_n$  of
$\CA$ by  
\begin{displaymath}   
\CA_n \, := \, \{\, a \in \CA \, : \, \forall k \geq n, [P_k,a] = 0 \, \}.
\end{displaymath}
Since each $s $ eventually commutes with the elements in the sequence
$(Q_n)$, the same is true for the sequence $(P_n)$ and we get that
\begin{displaymath} 
\forall s \in \CS \, \exists n_s \in \IN_0 \, :  \,  s \in
\CA_{n_s}
\end{displaymath}
In particular $\CS $ is contained in the union of the algebras
$\CA_n$, so this union is dense in $\CA$ and we only have to prove
that each $\CA_n$ is finite dimensional. Let us then fix a natural
number $n$. By the definition of $\CA_n$  we get a *-homomorphism
$\rho_n :  \CA_n \to B(P_nH)$ by $\rho_n(a) : = P_n a| P_nH = a|P_nH$,
 so in order
to prove that $\CA_n$ is finite dimensional it suffices to prove that
$\rho_n$ is faithful on $\CA_n$. Suppose now that $a$ in $\CA_n$
satisfies $\rho_n(a) = 0$, then since the separating vector $\xi $ is
in $P_nH$ we  get $0 = \rho_n(a) \xi = a \xi$  and $ a = 0$. The theorem
follows.
\end{proof}
\section{UHF C*-algebras}
A UHF C*-algebra is a special sort of AF C*-algebra where at each
stage the algebra $\CA_n$ is a full matrix algebra. In this section we
will then consider a UHF C*-algebra $\CA$ which is the norm closure of
an increasing sequence of finite dimensional full matrix algebras
$(\CA_n)_{n \in \IN_0} $ such that $\CA_0 = \IC I_{\CA}$. Since each
$\CA_n$ is a full matrix algebra there is an increasing sequence of
natural numbers $(m_n)_{n \in \IN}$ such that $\CA_n $ is isomorphic
to the full matrix algebra $\CM_{m_n}$. The assumption that the unit
of $\CA$ is the unit in all the algebras $\CA_n$ implies that there
must be natural numbers $d_n$ such that
\begin{displaymath}
\forall n \in \IN: \; \CA_n = \CA_{n-1} \otimes \CM_{d_n}
\end{displaymath}
In order to avoid trivial complications we will,  as in the previous
section, assume that the sequence $(m_n)_{n \in \IN} $ is strictly
increasing, i. e. all $d_n \geq 2$.
Based on this  we then have $d_1 = m_1$,  $m_n = d_1 \dots d_n$ and

\begin{displaymath} 
\CA_{n-1} = \CM_{d_1} \otimes \dots \otimes \CM_{d_{n-1}} \otimes \IC
I_{\CM_{d_n}} \subset   \CM_{d_1} \otimes \dots \otimes
\CM_{d_n} = \CA_{n}
\end{displaymath}

\noindent Going back to the notation from the section on AF C*-algebras  
we can now determine the space denoted $\CQ_n $ which is determined by
\begin{displaymath}
\forall n \in \IN: \; \CQ_n\; := \;  
\{ \pi_n(a) - \pi_{n-1}(a) \; : \; a \in  \CA \} 
\end{displaymath}
Since a full matrix algebra only has
one trace state, the normalized trace, it follows that 
 the  UHF C*-algebra $\CA$ also
only has one trace state, the restriction of which to $\CA_{n+1}  $ is
the tensor product of the trace states from   each of the factors
in the tensor product decomposition of $\CA_{n+1}$. Based on this, it
is for $n$ in $\IN_0$ 
 possible to describe $\CQ_{n+1} $ in terms of tensor products.
Let $\CM_{d_{n+1}}^{\circ}$ denote the elements in $\CM_{d_{n+1}}$ 
of trace zero,  then the product description of the trace state shows
that 
\begin{displaymath}
\forall n \in \IN_0: \; \CQ_{n+1}\; = \;  
\{ \pi_{n+1}(a) - \pi_n(a) \; : \; a \in  \CA \} =   \CM_{d_1} \otimes
\dots \otimes \CM_{d_n}\otimes \CM_{d_{n+1}}^{\circ}.
\end{displaymath}
In the proof of Theorem \ref{DAF} we introduced the constants
$c_{n+1}$. In the present case of a UHF algebra we  can 
compute the numbers $c_{n+1} $ exactly. Going back to the previous
section we find that $c_{n+1} $ is the maximal ratio between the
operator-norm and the 2-norm, with respect to the trace state,
 for elements in $ \CQ_{n+1}$. It is not
difficult to see that this maximum will be attained for an operator
which  is a rank one operator in 
$\CA_{n+1}$ provided that $\CQ_{n+1}$ contains a rank one operator. 
In order to see, that this is so,  one can consider an 
 operator  $f $ in
$\CA_{n+1}$ which is a product $ f = e_1 \otimes \dots \otimes e_{n+1}$ where
each $e_j$ is a matrix unit in $\CM_{d_j}$ and $\text{trace}(e_{n+1}) =
0 $. Such an operator $f$  must belong to 
$\CQ_{n+1}$ and be  a rank one operator in  $\CA_{n+1}$. Then for this
$f$ we get  
\begin{displaymath}
c_{n+1} \; = \; \frac{\|f\|}{\|f\|_2} \; = \frac{1}{(m_{n+1})^{-\frac{1}{2}}}
  \; = \; \sqrt{m_{n+1}}.
\end{displaymath}

\medskip

\noindent
We will now go back to the start Section \ref{AF} and take the
notation introduced there and apply it to the present
situation. Having the concrete values of the elements in the  sequence
$(c_n)_{n \in \IN} $ we can then modify  the proof of Theorem \ref{DAF}
 in accordance with the extra notation introduced just above,
 such that we can  obtain the following result.

\begin{Theorem} \label{DUHF}
Let $\CA$ be a UHF C*-algebra with an increasing sequence $(\CA_n)_{n
  \in \IN_0}$ of full matrix algebras such that $\CA_0 = \IC I_{\CA} $
and the union of the sequence is dense in $\CA$. With the notation
introduced  

\begin{itemize} 
\item[(i)] For an absolutely convergent 
 series $\sum \beta_n$ of non zero reals, the
  self adjoint  Dirac operator given by

\begin{displaymath}
D=\overset{\infty}{\underset{n=1}{\sum}}{(\beta_n)^{-1}\sqrt{m_n}Q_{n}}
\end{displaymath}  
will induce a metric on the state space for the weak*-topology.
This Fredholm module is 4-summable.
\item[(ii)] Given a positive $p < 2$ for $s >  \frac{2}{p} > 1$, the
  operator  $D$ given as  

\begin{displaymath}
D=\overset{\infty}{\underset{n=1}{\sum}}{(m_n)^sQ_{n}}
\end{displaymath}  
will induce a metric on the state space for the weak*-topology.
This Fredholm module is p-summable. 
\end{itemize}
\end{Theorem}
\begin{proof}
The first statement is a direct consequence of Theorem \ref{DAF}.
For the second statement one just has to remark that the assumption
$\CA_n \neq \CA_{n+1} $ implies that $m_{n+1} \geq 2m_n$, so $m_n \geq
2^n$. We can now do the estimation, first of the convergence of $\sum
\beta_k$
\begin{displaymath}
\sum_{k=0}^{\infty}\beta_{k+1}\, =\, \sum_{k=1}^{\infty} m_k^{(\frac{1}{2}
  - s)}\, \leq \,  \sum_{k=1}^{\infty} (2^{(\frac{1}{2}
  - s)})^k < \infty.
\end{displaymath}
The $p-$summability works by a similar sequence of estimates, one just
has to remember that the multiplicities of the eigenvalues
$(\text{dim}\CA_{n} - \text{dim}\CA_{n-1})$ is dominated by $(m_n)^2$.  
\begin{align*}
\text{tr}(I + D^2)^{-\frac{p}{2}} & \leq 1 + \sum_{n=1}^{\infty}
  (m_n)^{(2-ps)} \\
& \leq  \,1 +  \sum_{n=1}^{\infty} (2^{(2-ps)})^n \, < \,  \infty.
\end{align*}
\end{proof}
\bigskip

\bigskip
\section{On unbounded Fredholm modules for the Cantor set}
\noindent The usual Cantor set is
 a subset of the unit
interval in $\IR$ which is obtained by successive cuttings of
$2^{n-1}$ open
sub-intervals of length $(\frac{1}{3})^{n}$.
This space can also be considered as the compact
topological space $\Pi\IZ_2$ which is an infinite product of the
compact two
element group  $\{0,1\}$.
When viewed as an infinite product space it is quite easy to see that
it can equipped with  several 
inequivalent metrics, and it turns out that the Hausdorff dimension
for these metric spaces can attain any value in the interval $]0,
\infty[$.
  When we
want to emphasize that we are considering the classical Cantor set we
will call it the middle third Cantor set and denote it $\FC_{1/3}$,
the Hausdorff dimension of this space is $\frac{\log 2}{\log 3}$.
 To give the reader, who may
not be familiar with other representations of the Cantor set, an idea
of how to construct such a set, we
just mention that for a positive real number    $\gamma \in ]0,
\frac{1}{2}[$ it is possible to construct a homeomorphic copy of the
middle third Cantor set inside the unit interval by successive cuttings
of $2^{n-1}$  intervals of length $(1-2\gamma)\gamma^{n-1}$. This
space is denoted $\FC_{\gamma}$ and it's Hausdorff dimension is
$\frac{\log 2}{- \log{\gamma}}$.
\smallskip \noindent
For a C*-algebraist it is well known that the algebra of continuous
functions on the
Cantor set is an AF C*-algebra, but the C*-algebra carries only the
topological information, so all the geometry which comes from a
particular metric on the space has to be obtained from other sources.
The non commutative geometry as developed by Connes offers 
the tools to describe geometric data for C*-algebras in
general, and  Connes has tried to see how his theory can be applied to
the Cantor set $\FC_{1/3}$. This is done to some extent  in the book
\cite{Co3} ( IV 3.$\eps$),  and  in an unpublished  note \cite{Co4},
 we have had access
to,  Connes constructs an unbounded Fredholm module over the algebra 
$C(\FC_{1/3})$ which carries the
exact geometrical structure of the compact subset $\FC_{1/3}$ of $\IR$. 
\smallskip \noindent
We can not get so much information from the spectral triples we
associate to AF C*-algebras. It seems as the spectral triples we
 study for AF C*-algebras can reflect the structures which are
invariant under bi-Lipschitz mappings, and not much more.
What we do in this section is to  take the results from Section 2 and use the
extra structure to see what we can get for the Cantor C*-algebra.
This means that the representation of the algebra will be a fixed
standard representation and the eigenspaces of the Dirac operator will
be determined by the increasing sequence of finite dimensional
subalgebras. The spectral triples we can obtain this way are
 quite different
from the one Connes' constructs in \cite{Co3, Co4}, so we will just
shortly describe the major difference between Connes' module over 
$C(\FC_{1/3})$ and
the one we get from Section 2.

\subsection{Connes' spectral triple for the continuous functions on
  the Cantor set.}
\label{Connesspectral}

In his book  \cite{Co3}, ( IV 3.$\epsilon$ )  Connes describes a 
 Fredholm module over a general Cantor set  which is a subset of the
 interval $[0,1]$.  In
 the note \cite{Co4} Connes concentrates on the middle third Cantor
 set and constructs a spectral triple for the continuous functions on
 this set. The module is an infinite sum of one dimensional modules
 which are associated with the points in the unit interval which are
 either 0, 1 or an end point of some cut-out-interval. The Dirac
 operator has eigen values which 
reflects the distance between points which are either end
 points of a cut out interval or endpoints of an interval which is
 left  back 
 after a certain  number of cuts has been made.
This coding contains a
 surprising lot of the  geometrical data for  the space $\FC_{1/3}$.
In this construction, 
the module $H$ is a Hilbert space, which is an infinite  sum of one dimensional
modules over $C(\FC_{1/3})$. In particular $H$ has an orthonormal basis of
vectors, all  of which are joint eigenvectors for all the   elements
 in $C(\FC_{1/3})$.
In a language which makes sense for non commutative C*-algebras too, 
 one can say that Connes' module is a subrepresentation of
the reduced atomic representation  of $C(\CC)$. 

\smallskip \noindent 

\subsection{Some AF C*-algebra spectral triples
 for the continuous functions on
  the Cantor set.}
\label{AFspectral}

\noindent
The module, we have proposed as part of a spectral triple for  
 a general AF C*-algebras  with a
faithful trace state, yields for $C(\FC_{1/3})$  a
module which has no non trivial common eigenvectors but instead 
a separating and cyclic trace vector.

\smallskip
\noindent 
The special thing, for the  C*-algebra consisting of the continuous
functions on the Cantor set,  
is that it has a concrete description, and this  makes it possible
to perform detailed analysis. From Section 2 we know that
there is a lot of freedom in the choice of the eigenvalues for the
Dirac operator. This suggests that it would be interesting to see if it
is possible to describe the geometrical significance of some of the
possible choices of the eigenvalues.  We show below that
for any real $\gamma$ such that $0 < \gamma < 1$ we can use the
sequence $(\gamma^{-n+1})_{n \in \IN}$ as a sequence of eigenvalues, such
that the corresponding metric on the state space generates the
weak*-topology. Having this it seems natural to look for geometrical
consequences of this result and it turns out that if one equips the
usual Cantor set with the metric induced by a Dirac operator with the
sequence $(\gamma^{-n+1})_{n \in \IN}$ as  eigenvalues,
 then this compact  metric space will be of 
 Hausdorff dimension $\frac{\log 2}{-\log \gamma}$. We have
found a concrete realization of these fractals as a continuous family
of subsets of $\ell_1(\IN,\IR)$ and as subsets or $\IR^e$ for $e \in
\IN $ and $ e > \frac{\log2}{-\log \gamma}$.
All of this is of course based on the description of the topological
representation of the Cantor set as the compact Abelian group $\Pi
\IZ_2$. According to this description of the Cantor set we would like 
to mention that the discrete Abelian group  which is the dual of this
infinite product group is the infinite sum group which we denote 
$\oplus \IZ_2$, and the continuous functions on the Cantor set is just
the group C*-algebra for this discrete group. If one now 
equips this group with the length function
given by

\begin{displaymath}
\forall g = (g_n)_{n \in \IN} \in \oplus\IZ_2\,: \quad \ell(g) := \begin{cases}
\max\{ n \, : \, g_n = 1\} &  \rm{ if  }\, g \neq 0 \\
0 & \rm{ if  }\, g = 0, 
\end{cases}\end{displaymath}
\noindent
then the  spectral triple as constructed by Connes for the reduced
C*-algebra of a discrete group is - except for the size of the
eigenvalues - exactly the one we get from the AF-construction.

\smallskip
\noindent Before we state the theorem we would like to set up the frame, 
inside  which we will  work. First we will let $\CA$ denote the
C*-algebra C$(\Pi \IZ_2, \IC)$ and for $n$ in $\IN$ we 
 will let $e_n : \Pi \IZ_2 \to \IZ_2 $
 denote the coordinate mapping $ e_n((x_i)) = x_n$. Then $e_n$ is a
 self adjoint projection in $\CA$ and we will define a sequence of symmetries or self
 adjoint unitaries $(s_n)_{n \in \IN_0} $ in  $\CA$ by $s_0 = I$ and 
for $n >0,\;  s_n = 2e_n - I$. The sequence
 of finite dimensional subalgebras $(\CA_n)_{n \in \IN_0} $ is then defined
 such that $\CA_n$ is the least self adjoint complex algebra containing
the set $\{s_0, s_1, \dots , s_n\}$. Then $\CA_n$ is isomorphic to
$\IC^{2^n}$, the continuous functions on $2^n$ points. The union of
these algebras is denoted $\CA_{\infty}$ and is a self adjoint unital
subalgebra of C$(\Pi \IZ_2)$ which
 clearly separates the points in $\Pi \IZ_2$, so by
 Stone-Weierstrass' Theorem $\CA_{\infty} $ is dense in $\CA$ and
 $\CA$ is a unital AF C*-algebra. In order to apply Theorem \ref{DAF}
we have to fix a faithful state - or here a Borel probability measure
on $\Pi \IZ_2$ with support equal to $\Pi \IZ_2$. 
The natural choice is the measure which is determined in such a way
that the coordinate functions are treated symmetrically. This means
that the state, say $\tau $ is defined on $\CA_{\infty} $ by 
\begin{displaymath}
 \tau(I) = 1 \, \text{ and } \,  \tau(e_{i_1}e_{i_2}\dots e_{i_m}) =
 2^{- |\{i_1, i_2, \dots , i_m \}|}
\end{displaymath}
Then  $\tau$ is extended to all of $\CA$ by continuity and we will let
the corresponding regular Borel probability measure on
$\Pi \IZ_2$ be denoted by $\mu$. It should be noted that the
symmetry $s_n$ corresponds to the generator of the $n'$th summand in 
$\oplus \IZ_2$ and that this identification can be pursued to an
isomorphism which shows that the GNS-representation of $\CA$ on
$L^2(\CA, \tau)$ is nothing but the left regular representation of
C$^*_r(\oplus \IZ_2 )$.

\smallskip
\noindent
We can now formulate our theorem on Dirac operators for the Cantor set
using this notation and the definitions used in the set up for Theorem
\ref{DAF}
 
\begin{Theorem} \label{DCa}
Let $\CA$ denote the AF C*-algebra  {\rm C}$(\Pi \IZ_2)$, $ \CA_0 = \IC I_{\CA}, \,$ and for
$n \in \IN, \, \CA_n$ the subalgebra of $\CA$ 
generated by the $n$ first coordinate functions. 
\begin{itemize}
\item[(i)] If $(\alpha_n)_{n \in
  \IN_0}$ is a sequence of real numbers such that 
\begin{displaymath}
\alpha_0 = 0 \quad \text{ and } \sum_{n=1}^{\infty}\sup\{|\alpha_{n} -
\alpha_{i}|^{-1} \, : \, 0 \leq i \leq n-1 \} \, < \, \infty
\end{displaymath}
then the Dirac operator $D$ given by 
\begin{displaymath}
 D = \sum_{n=1}^{\infty} \alpha_n Q_n
\end{displaymath}
will generate a metric for the weak*-topology on the state space of 
{\rm C}$(\Pi \IZ_2)$
 and in particular a metric denoted 
$d_{(\alpha_n)}$ for the compact space $\Pi \IZ_2$
\item[(ii)] In the special case where there exists a real $\gamma$ such
  that $0 < \gamma < 1 $ and for $n \in \IN, \, \alpha_n = \gamma^{-n+1}$ 
  the conditions under (i) are fulfilled and the 
module will be $p$-summable for $p > \frac{\log 2}{- \log
    \gamma}$ and not for $p = \frac{\log 2}{- \log
    \gamma}$. The metric induced by $D$ on  $\Pi \IZ_2$ will in
  this case be denoted $d_{\gamma} $ and it will satisfy the following
  inequalities
\begin{align*}
\forall x \neq  y  \in \, \Pi \IZ_2 \text{ let } m(x,y)& = \min\{n \in
\IN \, : \, x(n) \neq y(n) \},  \text{ then } \\ 
2 \gamma^{m(x,y)-1} \leq d_{\gamma}(x,y) & \leq 2 \frac{\gamma^{m(x,y)-1}}{(1-\gamma)^2}
\end{align*}
\end{itemize}
\end{Theorem}
 \begin{proof}
In the proof of Theorem \ref{DAF} the important step is to get an
estimate of $\|\pi_{n+1}(a) - \pi_n(a)\| $ for an $a \in \CA$ such
that the   norm of the
commutator of $[D,a]$ is at most 1. Since the present situation is
much easier to deal with than the general one, we can get some rather
exact estimates of this type. So assume that we have an operator $a$
in $ \CD$, which means that $\|[D,a]\| \leq 1$. We start by examining
the expression $\|\pi_{n}(a) - \pi_{n-1}(a)\| $ in more details. First
we remark that -  since the algebra $\CA_n$ is represented faithfully on
the subspace $H_n$-  we have $\|\pi_{n}(a) - \pi_{n-1}(a)\| =
\|(\pi_{n}(a) - \pi_{n-1}(a))P_n\|$. Moreover $P_{n-1}(\pi_{n}(a) -
\pi_{n-1}(a))P_{n-1} = 0$ since    
\begin{align*} 
\forall b, c \in \CA_{n-1}: \, (\pi_n(a)b\xi,c\xi) \, &= \, (P_na\xi,
cb^{*}\xi) \\ &= (P_{n-1}a\xi, cb^{*}\xi) \, = \, (\pi_{n-1}(a)b\xi,c\xi).  
\end{align*}
A closer examination also shows, as
we shall see,  that
for this particular algebra also $Q_n(\pi_n(a)- \pi_{n-1}(a))Q_n = 0$.
This, last statement,
 follows from the facts that the space $F_n = Q_nH$ equals
$\CA_{n-1}s_n\xi$ and then for any $a$ in $\CA$ we have $\pi_n(a)-
\pi_{n-1}(a) \in \CA_{n-1}s_n$. Since $s_n^2 = I$ and $\CA $ is
commutative we get  $(\pi_n(a)- \pi_{n-1}(a))Q_nH  \subset P_{n-1}H$.
This all means that
\begin{displaymath}
P_n(\pi_n(a)- \pi_{n-1}(a))P_n = P_{n-1}(\pi_n(a)- \pi_{n-1}(a))Q_n + Q_n(\pi_n(a)- \pi_{n-1}(a))P_{n-1}. 
\end{displaymath}
For operators $\, x \,$ in $\CA$ which satisfy the relation 
\begin{displaymath}
x = P_{n-1}xQ_n + Q_nxP_{n-1} 
\end{displaymath}
we have $\|x\| = \max\{\|P_{n-1}xQ_n\|,\|Q_nxP_{n-1}\|\}$; this
follows from the C*-relation $\|x^*x\| = \|x\|^2$. We can therefore
obtain
\begin{align*}
\|\pi_{n}(a) - \pi_{n-1}(a)\| &=
\|(\pi_{n}(a) - \pi_{n-1}(a))P_n\|\\ 
 & = \max\{\|P_{n-1}(\pi_{n}(a) - \pi_{n-1}(a))Q_n\|,\|Q_n(\pi_{n}(a)
 - \pi_{n-1}(a))P_{n-1}\|\} \\
\end{align*}
\noindent
We are now going to relate the expression $\max\{\|P_{n-1}(\pi_{n}(a) - \pi_{n-1}(a))Q_n\|,\|Q_n(\pi_{n}(a)
 - \pi_{n-1}(a))P_{n-1}\|\}$ to the commutator  $[D,a]$ and its norm.
Since $D$
commutes with the projections $P_n$ and any operator in $\CA_n$
commutes with all the $Q_m$ for $m \geq n$ we get 
\begin{displaymath} \label{pnident}
\forall a \in \CD: \, P_n[D,a]P_n \, = \, [D, \pi_n(a)]P_n \, = \, P_n[D,
\pi_n(a)] \, = \, [D,\pi_n(a)]
\end{displaymath}
Based on this we get for an $a$ in $\CD$ and an $n$ in $\IN$
\begin{align*}
 P_{n-1}[D,a]Q_n  & = (DP_{n-1} - \alpha_nP_{n-1})aQ_n \\
&=(DP_{n-1} - \alpha_nP_{n-1})\pi_n(a)Q_n\\
&=(DP_{n-1} - \alpha_nP_{n-1})(\pi_n(a)-\pi_{n-1}(a))Q_n 
\end{align*}
and similarly 
\begin{align*}
 Q_n[D,a]P_{n-1}  & = Q_na(\alpha_nP_{n-1}-DP_{n-1}) \\
&=  Q_n\pi_n(a)(\alpha_nP_{n-1} - DP_{n-1}  )\\
&= Q_n(\pi_n(a)-\pi_{n-1}(a))( \alpha_nP_{n-1} - DP_{n-1}). 
\end{align*}
\noindent
The assumptions made on the sequence $(\alpha_n)_{n \in \IN_0}$ says
that the
series  $\sum_{n \in \IN} \beta_n$ of positive reals defined by $\beta_n
= \max\{|\alpha_n - \alpha_i|^{-1} \, : \, 0 \leq i \leq n-1 \}$ is
summable. When examining the expression $(\alpha_nP_{n-1}-DP_{n-1})$
one finds that the properties of the sequence  $(\alpha_n)$ implies 
 that this
operator is invertible on $P_{n-1}H$ with an inverse, say $B_n$, on
this space such that $\|B_{n}\| = \beta_n$. We can now make estimates
of $\max\{\|P_{n-1}(\pi_{n}(a) - \pi_{n-1}(a))Q_n\|,\|Q_n(\pi_{n}(a)
 - \pi_{n-1}(a))P_{n-1}\|\}$ by combining the previous computations.
Since it all depends on a max operation we may just as well assume
that $\|\pi_{n}(a) - \pi_{n-1}(a)\| = \|P_{n-1}(\pi_{n}(a) -
\pi_{n-1}(a))Q_n\|$ so
\begin{align} \label{beta} 
\|\pi_{n}(a) - \pi_{n-1}(a)\| &= \|P_{n-1}(\pi_{n}(a) -
\pi_{n-1}(a))Q_n\| \\
 &= \|B_n( DP_{n-1} - \alpha_nP_{n-1})(\pi_{n}(a) - \pi_{n-1}(a))Q_n\|
 \notag \\
&= \|B_nP_{n-1}[D,\pi_{n}(a) - \pi_{n-1}(a)]Q_n\|\notag \\
&= \|B_nP_{n-1}[D,a]Q_n\| \notag \\
& \leq \|B_n\|\|[D,a]\| \notag \\
& \leq \beta_n. \notag 
\end{align}

\medskip
\noindent
Now the arguments runs as in the proof of Theorem \ref{DAF}
and we get that the metric generated by the operator $D$ induces a
metric for the weak*-topology on the state space. For later use we
remark that 
\begin{equation} \label{stepest1}
\forall a \in \CD: \quad \|a - \tau(a)I\| \, = \, \|a - \pi_0(a)\| \,
\leq \, \sum_{i=1}^{\infty} \beta_i
\end{equation}
 
\begin{equation} \label{stepest2}
\forall a \in \CD \, \forall n \in \IN: \quad \|a - \pi_n(a)\| \,
\leq \, \sum_{i=n+1}^{\infty} \beta_i.
\end{equation}
 
\bigskip
\noindent We will now turn to the proof of item (ii), so we will
assume that for some $\gamma \in ]0,1[$ we have the eigenvalues of $D$
given by $ \alpha_0 =0$ and for $n \in \IN:  \alpha_n =
\gamma^{-n+1}$. In this setting the numbers $\beta_n $ are given by
\begin{displaymath}
\beta_n = (\gamma^{-n+1} - \gamma^{-n+2})^{-1} =
\frac{\gamma^{n-1}}{1-\gamma}, \, \text{ and } \sum\beta_n =(1-\gamma)^{-2}.
\end{displaymath}
  Having this we can for any  positive real $s$ 
compute the trace $\text{tr}(D^{-s})$, where as usual we disregard the
kernel of the selfadjoint positive operator $D$. The way the algebras
$\CA_n $ are defined shows that for $n \geq 1 $ we have  $\CA_n =
\CA_{n-1} \oplus (\CA_{n-1}s_n)$, 
hence   dim$Q_nH$ = dim$\CA_{n-1}s_n\, = $ dim$\CA_{n-1} = 2^{n-1}$
and we may compute the trace, 

\begin{equation} \label{trace} 
\text{tr}(D^{-s}) = \underset{k \to \infty}{\lim}\sum_{n=1}^k\gamma^{(n-1)s}2^{n-1}
= \underset{k \to \infty}{\lim}\frac{1 -
  (2\gamma^s)^{k}}{(1 -2\gamma^s)}
\end{equation}
\noindent
and it follows that this module is summable if and only if  $s > \log2
/ (-\log \gamma)$

\medskip
\noindent In order to determine the restriction to $\Pi \IZ_2$, 
of the metric which  $D$ induces on the set of regular Borel probability
measures, we first remark that for each $n
\in \IN$ we have $\pi_n(s_n) - \pi_{n-1}(s_n) = s_n$, so from
the computations done under item (i) we have 
\begin{align*}
[D,s_n]\, & =  \,  P_{n-1}[D,s_n]Q_n + Q_n[D,s_n]P_{n-1} \\
  & = (DP_{n-1} - \gamma^{-n+1}P_{n-1})s_nQ_n -  Q_ns_n(DP_{n-1} -
  \gamma^{-n+1}P_{n-1})
\end{align*}
\noindent
Since this expression again is of the form $Q_nxP_{n-1}+ P_{n-1}xQ_n$
the norm is the maximum of the norms of the summands. On the other
hand,
and as
described above,
 $s_nQ_n$ is a partial isometry with range projection $P_{n-1} $ and
 support projection $Q_n$. Hence 
 each of the summands above have norm equal to 
 $\|DP_{n-1} -
 \gamma^{-n+1}P_{n-1}\| = |\gamma^{-n+1}- 0| = \gamma^{-n+1}$ 
and we have got  
\begin{equation} \label{norm[D,s]}
\forall n \in \IN: \; \|[D,s_n]\| =  \gamma^{-n+1}.
\end{equation}
which means that 
\begin{equation} \label{SninD}
\forall n \in \IN: \;  \gamma^{n-1} s_n \in \CD.
\end{equation}
We will now determine the quantitative  
consequence of (\ref{SninD}) for the metric
$d_{\gamma}$ For  $x \neq  y \in \Pi \IZ_2$ we recall that
$m(x,y)$ is the least natural number for which the coordinates satisfy
$x(n) \neq y(n)$. Let then $\chi_x$ and $ \chi_y$ denote the states on
$\CA$ or probability measures on $\Pi \IZ_2$
 which are the point evaluations at  $x$
and $y$.  We then get 
\begin{equation} \label{dgeq}
d_{\gamma}(x,y) \, \geq \, |\chi_x(\gamma^{m(x,y)-1}s_{m(x,y)}) -
\chi_y(\gamma^{m(x,y)-1}s_{m(x,y)})| = 2\gamma^{m(x,y)-1}. 
\end{equation}

\smallskip
\noindent In order to obtain an inequality in the opposite direction
we have to return to the series of  inequalities (\ref{stepest2}). Given
$x \neq  y$  in $\Pi \IZ_2$,  an $\eps > 0$ 
and an $a$ in $\CD$ such that $|\chi_x(a) -
\chi_y(a)| \geq d_{\gamma}(x,y) - \eps$, then 
we get for $n \in \IN_0 \text{ and } n < m(x,y)$ that 
$\chi_x(\pi_n(a)) = \chi_y(\pi_n(a))$, so by (\ref{stepest2})    
\begin{align} \label{dleq}
d_{\gamma}(x,y) \, &\leq \, |\chi_x(a) -
\chi_y(a)| + \eps  \\
& = |\chi_x(a - \pi_{(m(x,y)-1)}(a)) - \chi_y(a-\pi_{(m(x,y)-1)}(a))|
+ \eps  \notag \\
& \leq 2 \|a - \pi_{(m(x,y)-1)}(a)\| + \eps \, \text{ and since }
\beta_i = \frac{\gamma^{i-1}}{1-\gamma}  \notag \\
& \leq  2 \sum_{i=m(x,y)}^{\infty}\frac{\gamma^{i-1}}{(1 - \gamma)}
 + \eps \notag \\
& =2\frac{\gamma^{m(x,y)-1}}{(1-\gamma)^2} + \eps \notag. 
\end{align}
We can then conclude that 
\begin{equation} \label{dcomp}
\forall x, y \in \Pi \IZ_2: \quad 2\gamma^{m(x,y)-1} \,  \leq \, 
d_{\gamma}(x,y) \, \leq \, 2\frac{\gamma^{m(x,y)-1}}{(1-\gamma)^2}
\end{equation}
and the theorem follows.
\end{proof}

\bigskip
\noindent
The natural question raised by Theorem \ref{DCa},  is  whether
there exist {\em ``natural''}  geometrical representations of the
Cantor sets
corresponding to all the values of $\gamma \in ]0,1[$. At a first
sight it is not even clear what the question means, but for the type
of metric spaces  we consider the relevant equivalence concept seems
to be that of {\em bi-Lipschitz }. We will recall this concept
\cite{Fa} in the definition below.

\begin{Definition}
Let $( S, \mu)  $ and $(T, \nu) $ denote two metric spaces and $F : S
\to T$ a map. The map $F$ is said to be bi-Lipschitz if it is  bijective
and  there
exist strictly positive constants $k, K$ such that
\begin{displaymath}
\forall x, y \in S: \, k\mu(x,y)\, \leq \, 
\nu(F(x),F(y)) \, \leq \,  K\mu(x,y).
\end{displaymath}
If $S = T$ and the identity mapping on $S$ is bi-Lipschitz 
then we say that $\mu$ and $\nu$ are equivalent metrics on $S$.
\end{Definition}

\noindent
Having this definition we will first remark that it is quite clear
from the last inequalities in the Theorem \ref{DCa} that for a pair of
different values of $\gamma$ the corresponding metrics are
inequivalent. A Closer look at  the mentioned inequalities makes it
apparent that for a given $\gamma \in ]0,1[ $ there may be a
representative for the metrics in the equivalence class containing
$d_{\gamma}$ which may look more attractive than the rest. The metric
which we have chosen to like the most is given by 

\begin{Definition}
Let  $\Pi \IZ_2$ denote the Cantor set. For any $\gamma \in
]0,1[$ we let $\delta_{\gamma} $ denote the metric on
$\Pi \IZ_2$
given by
\begin{displaymath}
\forall x, y \in \Pi \IZ_2: \; \delta_{\gamma}(x,y) =  
\sum_{n=1}^{\infty} |x(n) - y(n)| \gamma^{n-1} (1-  \gamma)
\end{displaymath}
\end{Definition}

\noindent
We will now compute the bounds for the equivalence between
$\delta_{\gamma} \text{ and } d_{\gamma}$.

\begin{Proposition} \label{compare}
For any $\gamma \in ]0,1[$ and any $x, y \in \Pi \IZ_2$
\begin{itemize}
\item[(i)]$ \;
 \gamma^{(m(x,y)-1)}(1-\gamma) \, \leq \, \delta_{\gamma}(x,y) \, \leq
 \,  \gamma^{(m(x,y)-1)}$
\item[(ii)] $ \; 
2 \delta_{\gamma}(x,y) \, \leq \,  d_{\gamma}(x,y) \, \leq \,
\frac{2}{ (1-\gamma)^3}\delta_{\gamma}(x,y).$
\item[(iii)]
For any two different values of $\gamma$ the corresponding metrics are
inequivalent.
\end{itemize}
\end{Proposition}

\begin{proof}
The sum formula for geometric series yields (i) right away. 

\noindent 
For (ii) we see that a
 combination of (i) with the results of  Theorem \ref{DCa} (ii)
 yields the
equivalence as stated. The very same reference shows that different
$\gamma$'s give inequivalent metrics.
\end{proof}

\subsection{Continuous families of Cantor sets of different dimensions}
\label{contCa}

Based on the construction of the metrics $\delta_{\gamma}$ it is clear
that the metric space $(\Pi \IZ_2, \delta_{\gamma}) $
has an {\em isometric embedding } onto a compact subset  of the
Banach space $\ell^1(\IN,\IR)$ of absolutely summable sequences of
real numbers. It is quite easy to compute the dimension of this space
as $\log 2 / (- \log \gamma)$, but it  was less obvious to us that for
the natural number $d$ defined as $d = \lfloor \log 2 / (- \log
\gamma) \rfloor +1 $, it is possible to construct a bi-Lipschitz map of
$(\Pi \IZ_2, \delta_{\gamma}) $ onto a compact subset of
$\IR^d$ equipped with the usual metric.
 If $0 <\gamma < 1/2 $, the image is the usual Cantor set obtained by
cutting middle intervals of the unit interval, such that at step $n$
we cut $2^{n-1} $ intervals each of length $(1-2\gamma)\gamma^{n-1}$. 
If $d > 1$ then the image is a product in $\IR^d$ of 
$d$ copies of the Cantor set 
where the cutting lengths on the sides are
$(1-2\gamma^d)\gamma^{d(n-1)}$.
 
\smallskip
\noindent
We will now define  the embedding of
$(\Pi \IZ_2, \delta_{\gamma}) $ into $\ell^1(\IN, \IR)$.
\begin{Definition}
Let $\gamma \in ]0,1[$ then $f_{\gamma}: (\Pi \IZ_2,
\delta_{\gamma}) \to   \ell^1(\IN,\IR)$ is defined by 
\begin{displaymath}
\forall x \in \Pi \IZ_2 \, \forall n  \in \IN:  \;
  f_{\gamma}(x)(n) := \gamma^{n-1}(1-\gamma)x(n).
\end{displaymath}
The image in $\ell^1$ of the metric set $(\Pi \IZ_2, \delta_{\gamma})$ by
$f_{\gamma} $  is denoted $\CC_{\gamma}$
\end{Definition}

\begin{Proposition}
For any $\gamma \in ]0,1[$ the map $f_{\gamma}$ 
 is an isometry of $(\Pi \IZ_2,
\delta_{\gamma})$ onto $(\CC_{\gamma}, \|.\|_1)$
\end{Proposition}
\begin{proof}
\begin{displaymath}
\forall x,y \in \Pi \IZ_2:
\|f_{\gamma}(x)-f_{\gamma}(y)\|_1 =
\sum_{n=1}^{\infty}|x(n)-y(n)|\gamma^{n-1}(1-\gamma) =
\delta_{\gamma}(x,y)
\end{displaymath}
\end{proof}

\begin{Remark}
We have chosen to let the embedding take place in the space,
$\ell^1(\IN, \IR)$ because of the isometric embedding. 
On the other hand - as sets - we have      $\ell^1(\IN, \IR) \subset
\ell^p(\IN, \IR)
 \subset \ell^{\infty}(\IN, \IR)$ for any $p \in ]1, \infty[$, so  it
 might be reasonable to look at the other embeddings too. It turns out
 that the metrics coming from the other norms are all equivalent on
 the set $\CC_{\gamma}$. An elementary computation shows these equivalences 
\begin{align*}
\forall x, y \in \Pi \IZ_2\, \forall p \geq 1: \; \\
\gamma^{(m(x,y)-1)}(1-\gamma)  & = \|f_{\gamma}(x) -
f_{\gamma}(y)\|_{\infty} \\ & \leq  \|f_{\gamma}(x) -
f_{\gamma}(y)\|_p \\ & \leq  \|f_{\gamma}(x) -
f_{\gamma}(y)\|_1  \\ &\leq \gamma^{(m(x,y)-1)}.
\end{align*}
\end{Remark}

\noindent
The Hausdorff dimension of the metric spaces $\CC_{\gamma}$ can be
computed in the same way as it is done for the usual Cantor set
contained in the unit interval. We will base our computations on the
one K. Falconer offers in \cite{Fa} on the pages 31 - 32. We will use
the same notation as Falconer, so  $\mathscr{H}^s(\CC_{\gamma}) $
denotes the $s$ dimensional Hausdorff measure of $\CC_{\gamma}. $ 
The Hausdorff
dimension, say $t$, is characterized by the fact that
$\mathscr{H}^s(\CC_{\gamma}) = 0,  $ 
if $s > t$ and the value is
infinite if $s < t$.
 Let us then fix a $\gamma \in ]0,1[$ and define $t:= \frac{\log 2}{-
  \log \gamma}$. We will then show that the $t$ dimensional 
 Hausdorff measure of $\CC_{\gamma}$ is a strictly positive real
 number, 
so the
 Hausdorff dimension of $\CC_{\gamma} $ is $ \frac{\log 2}{-
  \log \gamma}$.      

\smallskip
\noindent In order to ease the translation of the proof from \cite{Fa}
we will give a short description of the {\em standard intervals } in
$\CC_{\gamma} $ and describe some of their properties in a lemma.  For
$n \in \IN$ we will let $\CS_n$ denote all the points in $\Pi \IZ_2$
whose coordinates are all zero from coordinate number 
$n+1$ and onwards,  or formally
$\CS_n = \oplus_{i=1}^n \IZ_2 \subset \Pi \IZ_2$. The union of the
sets $\CS_n$ is denoted $\CS$. It seems to be convenient to introduce
the projection mappings $\pi_n : \Pi \IZ_2 \to \CS_n$ which for an $x
\in \Pi \IZ_2$ replaces all the coordinates of $x$ from the number
$n+1 $ and onwards by $0$.  For an $s$ in $\CS$ and an $n$ in $\IN$ we
then define the standard interval $V(s,n) $ by
\begin{displaymath}
\forall s \in \CS \, \forall n \in \IN: \quad V(s,n) = 
\begin{cases} \emptyset, & \text{ if } s \notin \CS_n \\
\{x \in \Pi \IZ_2\, : \, \pi_n(x) = s \}, & \text{ if } s \in \CS_n 
\end{cases}
\end{displaymath}
\noindent
This description of $V(s,n)$ works independently of the chosen
$\gamma$, but in the following lemma, where we list some properties of
the standard intervals, we also state some properties which relate to
the present metric, so we will now consider the intervals as subsets
of  $\CC_{\gamma}$.

\begin{Lemma} \label{standlem}

The standard intervals $V(s,n)$ of $\CC_{\gamma}$ have the following
properties:

\begin{itemize}
\item[(i)] $\forall n \in \IN$ there are exactly $2^n$ nonempty
  standard intervals of the form $V(s,n)$,  and they are indexed by
  the points in $\CS_n$ by $\{V(s,n)\, : \, s \in \CS_n \}$. Further
  the sets are pairwise disjoint and their union  equals $\Pi \IZ_2$. 
  \item[(ii)] The $2^n$ sets $\{V(s,n)\, : \, s \in \CS_n \}$ are open
  and closed and  each one has diameter equal to
  $\gamma^n$ in $\CC_{\gamma}$.
  
\item[(iii)] For any subset $U $ of $\CC_{\gamma}$ of diameter $|U| <
  (1-\gamma)$ there exists a standard interval $V(s,n)$ which contains
  $U$ and such that it's diameter satisfies $|V(s,n)| \leq
  |U|/(1-\gamma)$.
  
\item[(iv)] For  $t = \frac{\log 2}{- \log \gamma}$, for any 
$n, m \in \IN $ and for any $s_0$ in $ \CS_n$
\begin{displaymath}
2^{-n} \, = \, |V(s_0,n)|^t \, = \, \sum_{s \in \{x \in \CS_{m+n} \, : \, \pi_n(x) =
  s_0 \}} |V(s,n+m)|^t.
\end{displaymath}
\end{itemize}
\end{Lemma}
\begin{proof}
The content of (i) follows from the fact that $\CS_n$ has exactly
$2^n$ points. With respect to (ii) the topological content is
obvious since the coordinate mappings map into a two-point space. The
diameter estimate follows from the inequality below
\begin{displaymath}
\forall n \in \IN \, \forall s \in \CS_n \, \forall x, y \in V(s,n):
\; \delta_{\gamma} (x,y) \,=\, \sum_{i = n+1}^{\infty} |x_i -
y_i|\gamma^{i-1}(1-\gamma)  \, \leq \, \gamma^n.
\end{displaymath}
It is evident that for any $x$ in $V(s,n)$ there is a $y$ in $V(s,n)$
such that $\delta_{\gamma}(x,y) = \gamma^n$, so the maximum distance
is not only attained, but it can be attained from  any point in the
set !
The result of (iii) is not so obvious and demands some more
computations. Let then $U$ be given such that $|U| < 1- \gamma$. We
will fix a point $u$ in $U$, then for any $x$ in $U$ we have 
\begin{displaymath}
|U| \geq \delta_{\gamma}(u,x) = \sum_{i = m(u,x)}^{\infty}|u_i - x_i|
\gamma^{i-1}(1-\gamma) \geq \gamma^{(m(u,x)-1)}(1-\gamma).
\end{displaymath}
This shows that
\begin{equation*}
\forall x \in U: \; m(u,x) - 1 \geq \lceil
\log(\frac{|U|}{1-\gamma})\frac{1}{\log \gamma} \rceil
\end{equation*}
Hence for 
\begin{equation*}
n =  \lceil
\log(\frac{|U|}{1-\gamma})\frac{1}{\log \gamma} \rceil  \text{ and }
 s = \pi_n(u): \; U \subset  V(s,n)
\end{equation*} 
and since the diameter of $V(s,n) = \gamma^n$  we get that 
\begin{displaymath}
|V(s,n)| = \gamma^{\lceil
\log(\frac{|U|}{1-\gamma})\frac{1}{\log \gamma} \rceil} 
\leq \gamma^{\log(\frac{|U|}{1-\gamma})\frac{1}{\log \gamma}} =
\frac{|U|}{1-\gamma},
 \end{displaymath}
so (iii) follows.
With respect to (iv) we  remark that $t $ is defined such that
$\gamma^t = \frac{1}{2}$. In order to prove (iv) it is by the
induction principle enough to consider the case where m = 1. So let
$n$ in $\IN$ and $s$ in $\CS_n$ be given then $|V(s,n)|^t =
(\gamma^n)^t =  2^{-n} $.
Each such interval  is divided
into two intervals $V(s_0, n+1) \text{ and } V(s_1, n+1) $ 
 at the level $n+1$, and by the computations just performed we have
\begin{displaymath}
|V(s_0, n+1)|^t +  |V(s_1, n+1)|^t = 2^{-n-1} + 2^{-n-1} = 2^{-n} =
|V(s,n)|^t
\end{displaymath}
The lemma follows.
\end{proof}

The result on the Hausdorff dimension for $\CC_{\gamma}$ can  now be proved.

\begin{Theorem}
For any $\gamma$ in $]0,1[$, and for $t = \frac{\log 2}{-\log \gamma}$
the Hausdorff measure of the compact metric space $\CC_{\gamma}$
satisfies 
$(1-\gamma)^t  \leq \mathscr{H}^t(\CC_{\gamma}) \leq 1$.
\end{Theorem}

\begin{proof}
With respect to the upper estimate, we take an   $n$ in $\IN$ and
consider  the collection of
intervals $\{V(s,n) \, : \, s \in \CS_n \}$. This  constitutes a covering of
$\CC_{\gamma} $ by $2^n$ sets of diameter $\gamma^n$. By Lemma
\ref{standlem} (iv) or by the definition of $t$, 
we have $\sum_{s \in \CS_n}|V(s,n)|^t =1$, so 
$\mathscr{H}^t(\CC_{\gamma}) \leq 1 $. 

\noindent
The inequality the other way is determined after a couple of reduction
steps. As Falconer mentions one can without loss of generality
restrict the attention to covering families $(U_i)$, where each set
$U_i$ is an open set. Since $\CC_{\gamma} $ is compact this means, in
turn, that we can restrict to finite families $(U_i)$. Let then $0 <
\delta < 1 - \gamma $ and  a finite $\delta$-cover $(U_i)$  be given.
For each $i$ we can then by Lemma \ref{standlem} (iii) 
find a natural number $n_i$ and  an  $s_i$ in
$\CS_{n_i}$ such that $U_i \subset V(s_i, n_i)$ and $(1-\gamma)|V(s_i, n_i)|
\leq |U_i|$. Further since there are only finitely many $n_i$ we can
define $m = \max\{(n_i)\}$ and by Lemma \ref{standlem} (iv)
 we can start estimating  
\begin{align*}
\sum |U_i|^t &\geq (1-\gamma)^t \sum |V(s_i, n_i)|^t   \\
& = (1-\gamma)^t\sum_{i} \sum_{s \in \CS_m \text{ and } \pi_{n_i}(s) =
  s_i}|V(s,m)|^t \\
& \geq  (1-\gamma)^t \sum_{s \in \CS_m}|V(s,m)|^t \\
& = (1-\gamma)^t, 
\end{align*} 
so 
 $\mathscr{H}^t(\CC_{\gamma}) \geq  (1-\gamma)^t$
and the theorem follows.
\end{proof}
\noindent
For each $\gamma$ in $]0,1[$ we will find a suitable natural number
$e_{\gamma}$ such that we can embed $\CC_{\gamma} $  into  
$\IR^{e_{\gamma}}$  via a bi-Lipschitz mapping. 
Since bi-Lipschitz maps preserve Hausdorff dimensions it follows that we
can obtain finite dimensional representations of all the
$\CC_{\gamma}$'s. 

\begin{Theorem}
Let $\gamma $ be in $]0,1[$,  $e_{\gamma} := \lfloor \frac{\log 2}{-\log
  \gamma} \rfloor +1 $ and let $F_{\gamma} : \CC_{\gamma} \to
\IR^{e_{\gamma}}$ be defined via  it's $i$'th coordinate
function,  $F^i_{\gamma}$, as 
\begin{displaymath}
\forall i \in \{1, \dots , e_{\gamma}\}\, \forall x  \in \Pi \IZ_2 :
 \, F^i_{\gamma}(x) = \sum_{p=1}^{\infty}x(i + (p-1)e_{\gamma})(\gamma^{
 e_{\gamma}})^{p-1}(1-\gamma^{e_{\gamma}}).
\end{displaymath}
Then $F_{\gamma}$ is a bi-Lipschitz continuous mapping of $\CC_{\gamma}
$ onto it's image as subset of $\IR^{e_{\gamma}}$.

\noindent
If $\gamma < 1/2$ then $e_{\gamma} =1 $ and 
the image is the usual $\gamma$-Cantor subset of the
unit interval which can be obtained by successive cuttings of $2^{n-1}$ middle
intervals of length  $\gamma^{n-1}(1 -2 \gamma)$. If $\gamma \geq \frac{1}{2}$
then the image is the product of $e_{\gamma} $ copies of the one
dimensional $\gamma^{e_{\gamma}}$-Cantor subset
$\FC_{\gamma^{e_{\gamma}}}$
 of the unit interval. 
  \end{Theorem}

\begin{proof}
It is obvious that $F_{\gamma} $ is continuous, but the bi-Lipschitz
property can only be seen after a few computations, but before we
start the computations we mention that the symbol $\|z\|$
means the usual Euclidian norm in $\IR^{e_{\gamma}}$. Other norms will
be indicated by a subscript. The
continuity and the Lipschitz property for $F_{\gamma}$ follow from the
following inequalities.

\begin{align*}
 & \forall x, y \in \CC_{\gamma}: \|F_{\gamma}(x) - F_{\gamma}(y)\| \\ 
& \leq \, \|F_{\gamma}(x) - F_{\gamma}(y)\|_1 \notag \\
&  \leq \, \sum_{i=1}^{e_{\gamma}} \sum_{p=1}^{\infty} |x(i +
 (p-1)e_{\gamma}) - y(i + (p-1)e_{\gamma})|(\gamma^{
 e_{\gamma}})^{p-1}(1-\gamma^{e_{\gamma}}) \notag \\
& = \sum_{i=1}^{e_{\gamma}} \sum_{p=1}^{\infty} \frac{\gamma^{1-i}(1-
  \gamma^{e_{\gamma}})}{1-\gamma}|x(i +
 (p-1)e_{\gamma}) - y(i + (p-1)e_{\gamma})|\gamma^{(i +
   (p-1)e_{\gamma} -1)}(1-\gamma) \notag \\
& \leq \sum_{i=1}^{e_{\gamma}} \sum_{p=1}^{\infty}
\frac{\gamma^{(1-e_{\gamma})}
}{1-\gamma}|x(i +
 (p-1)e_{\gamma}) - y(i + (p-1)e_{\gamma})|\gamma^{(i +
   (p-1)e_{\gamma} -1)}(1-\gamma) \notag \\
& = \frac{\gamma^{(1-e_{\gamma})}
}{1-\gamma}\delta_{\gamma}(x,y) \notag
\end{align*}

\noindent Before we start to prove that $F_{\gamma}$ has a Lipschitz
inverse we remind the reader that the choice of $e_{\gamma}$  as 
consequence has, that $\gamma^{e_{\gamma}} <  \gamma^{(\frac{\log 2}{ - \log
    \gamma })} = \frac{1}{2}$.
Let us fix a pair $x, y$ in $\CC_{\gamma}$ and let 
$i $ in $\{1, \dots \ e_{\gamma} \}$ be chosen such that $ i
\equiv m(x,y) $ mod $e_{\gamma}$ then 
we will first determine $p$ such that $i + (p-1)e_{\gamma} = m(x,y)$.
Since $i = e_{\gamma} $ is possible, it turns out that 
\begin{displaymath}
p= \lceil \frac{m(x,y)}{e_{\gamma}}\rceil.
\end{displaymath}
We will base the estimates below on the $i$'th coordinate function,
still the same $i$.
We will start the estimation on $|F^i_{\gamma}(x)-
F^i_{\gamma}(y)|$, by using the knowledge on $p$ just obtained
\begin{align*}
 |F^i_{\gamma}(x)- F^i_{\gamma}(y)| & = | \sum_{p = \lceil
   \frac{m(x,y)}{e_{\gamma}}\rceil}^{\infty}( x(i +
 (p-1)e_{\gamma}) - y(i + (p-1)e_{\gamma}))(\gamma^{
 e_{\gamma}})^{(p-1)}(1-\gamma^{e_{\gamma}})|  \\
& \geq  (\gamma^{e_{\gamma}})^{(\lceil
   \frac{m(x,y)}{e_{\gamma}}\rceil -1)}(1-\gamma^{e_{\gamma}})
 -\sum_{p = \lceil 
   \frac{m(x,y)}{e_{\gamma}}\rceil +1}^{\infty}(\gamma^{
 e_{\gamma}})^{(p-1)}(1-\gamma^{e_{\gamma}}) \notag \\
  & =  (\gamma^{e_{\gamma}})^{(\lceil
   \frac{m(x,y)}{e_{\gamma}}\rceil -1)}(1-\gamma^{e_{\gamma}}) -
 (\gamma^{e_{\gamma}})^{\lceil \frac{m(x,y)}{e_{\gamma}}\rceil} \notag
 \\
& = (\gamma^{e_{\gamma}})^{(\lceil
   \frac{m(x,y)}{e_{\gamma}}\rceil -1)}(1-2\gamma^{e_{\gamma}}) \notag \\
& = (\gamma^{e_{\gamma}})^{(p-1)}(1-2\gamma^{e_{\gamma}}) \notag \\
& = (\gamma^{(m(x,y)-i)}(1-2\gamma^{e_{\gamma}}) \notag \\ 
& \geq \gamma^{m(x,y)}(1- 2\gamma^{e_{\gamma}}) \notag
\end{align*}
The estimates on $F_{\gamma}$ can now be performed using Proposition
\ref{compare} 
\begin{align*} 
\|F_{\gamma}(x) - F_{\gamma}(y)\| \, 
& \geq \, |F^i_{\gamma}(x) - F^i_{\gamma}(y)| \\
& \geq \,  \gamma^{m(x,y)}(1- 2\gamma^{e_{\gamma}})
\notag \\
& \geq \, (1-2\gamma^{e_{\gamma}}) \gamma \delta_{\gamma}(x,y) \notag
\end{align*}

We have now proved that each map 
 $F_{\gamma} $ is  bi-Lipschitz  and
we only have to study the image $F_{\gamma}(\CC_{\gamma})$. We will
stick to the case where $\gamma = 1/3$ and consider an $x$ in $\Pi \IZ_2$
\begin{align*}
F_{\gamma}(x) &  = \sum_{n=1}^{\infty} x(n)(\frac{1}{3})^{n-1}(1 -
  \frac{1}{3}) \\
& = \sum_{n=1}^{\infty} x(n) 2(\frac{1}{3})^n.
\end{align*}
Since $ x(n) $ is in the set $\{0,1\}$ we see that the image of
$F_{\gamma}$ is exactly the set of points in the unit interval which
in the base 3, can
written using only  digits from the set $\{0, 2\}$. This is a well
known characterization of the standard Cantor set, $\FC_{1/3}$.
The theorem follows.
\end{proof}

\subsection{ On the Gromov--Hausdorff distance between $\CC_{\gamma}
  \text{ and } \CC_{\mu}$. } 

\noindent
Hausdorff has defined  a metric on the
closed subsets of a compact metric space and in this way obtained  a new
compact metric space. In the book \cite{Gr} M. Gromov extends this
idea and defines a distance between any pair of compact metric
spaces. This metric on compact metric spaces is denoted the {\em
  Gromov--Hausdorff distance. } 
In order to  provide a frame, inside which it makes sense
to speak about convergence of a family of 
non commutative compact metric spaces,  Marc
A. Rieffel has studied many aspects of 
the Gromov--Hausdorff distance between non
commutative compact metric spaces, and published a memoir in the series of
the American Mathematical Society on these matters \cite{Ri6}. 
After the  first presentation of our  application of the C*-algebraic
approach to the Cantor set, Marc A. 
Rieffel asked us, if we could compute the Gromov--Hausdorff
  distance between the different metric spaces
$\CC_{\gamma} $ obtained from the Cantor set. We can not answer this
question completely but below we can give an upper bound on the
distance between  two such spaces  $\CC_{\gamma}
  \text{ and } \CC_{\mu}$. The estimate is based on the 
isometric embedding of each $\CC_{\gamma } $
in $\ell^1(\IN) $  as constructed in the previous subsection.

\noindent
 In order to make
this computation we first  recall the definition of the Hausdorff
distance between closed subsets of a metric space and then the
definition  of the Gromov--Hausdorff
distance between metric spaces. 
\begin{Definition}
For closed subsets $X \text{ and } Y $ of a compact 
metric space $(Z, d)$  the Hausdorff distance between $X$ and $Y$
is given by

\begin{displaymath}
\text{dist}^d_H(X,Y) = \inf\{ r>0 \, :\, \forall x \in X \exists y \in
Y\,  d(x,y) \leq r \text{ and }  \forall y \in Y \exists x \in
X\,  d(y,x) \leq r \, \} 
\end{displaymath}
\smallskip
\noindent
For compact metric spaces $(X,d_x)$ and $(Y,d_y)$ let
$X \Dot{\cup} Y$ denote the disjoint union of the sets and let
$\CM(d_x, d_y) $ denote the set of all metrics on this space such that
the restriction of any of these metrics to each of the subsets $X$ and
$Y$ agrees with the given metric on that space. The Gromov--Hausdorff
distance between $ X$ and $Y$ is then given by
\begin{displaymath}
\text{dist}_{GH}((X,d_x),(Y,d_y)) = \inf\{ \text{dist}^d_H(X,Y ) \, : \, d \in
\CM(d_x, d_y)\, \}
\end{displaymath}
\end{Definition}
\begin{Proposition}
Let $0 < \mu < \gamma < 1 $ then
\begin{displaymath} 
\rm{dist}_{GH}(\CC_{\gamma},\CC_{\mu} )\, \leq \,2 \frac{\gamma  -
  \mu }{1 - \gamma}.
\end{displaymath}
\end{Proposition}
\begin{proof}
The spaces $\CC_{\gamma} $ are all subsets of $\ell^1(\IN)$, but they
are not disjoint inside this metric space since for instance $0$ belongs
to all of them. This deficiency can easily be repaired by considering
the normed  space $E = \IR \times \ell^1(\IN, \IR)$ which is equipped
with the norm $\|(t,x)\|:= \max\{|t|, \|x\|_1\}$ and by the
construction  of isometric copies $\CD_{\gamma} $ of $\CC_{\gamma}$ in
$E$. This is done by defining
  
\begin{displaymath}
\CD_{\gamma} := \{( \gamma, (x_n\gamma^{n-1}(1-\gamma))_{n \in \IN})
\in E
\, : \, x_n \in \{0,1\}\, \}
\end{displaymath}
\noindent
We first remark that since $  \mu < \gamma $ we must have
$\gamma^{n-1}(1-\gamma) \leq \mu^{n-1}(1-\mu)$ for at least $n=1$ and
we define $m$ as the largest natural number such that the inequality
above is satisfied. Then $\gamma^{n-1}(1-\gamma) > \mu^{n-1}(1-\mu)$
for $ n > m$. We can now make  an  estimate of the Hausdorff
distance between $\CD_{\mu} \text{ and } \CD_{\gamma}$. Let $u = (\mu,
(x_n\mu^{n-1}(1-\mu))_{n \in \IN})$ 
be a point in $\CD_{\mu}$ then for $v = (\gamma,
(x_n\gamma^{n-1}(1-\gamma))_{n \in \IN})$ in $\CD_{\gamma}$, with the
same sequence $(x_n)$   we can perform the estimates

\begin{align*}
&\|u-v\|_E \\ & = \max\{ \gamma - \mu, \sum_{n =1 }^{\infty} x_n |
\gamma^{n-1}(1-\gamma) - \mu^{n-1}(1-\mu)| \}  \\
&= \max\{ \gamma - \mu, \sum_{n = 1 }^{m} x_n( 
 \mu^{n-1}(1-\mu)- \gamma^{n-1}(1-\gamma)) + \sum_{n = m+1 }^{\infty} x_n( 
\gamma^{n-1}(1-\gamma) - \mu^{n-1}(1-\mu)) \} \\
&\leq  \max\{ \gamma - \mu, \sum_{n = 1 }^{m}  
 \mu^{n-1}(1-\mu)- \gamma^{n-1}(1-\gamma) + \sum_{n = m+1 }^{\infty}  
\gamma^{n-1}(1-\gamma) - \mu^{n-1}(1-\mu) \} \\
&= \max\{ \gamma - \mu, 2( \gamma^m - \mu^m)\}
\end{align*}

\noindent
By definition of $m$ we have $\gamma^{m-1}(1- \gamma) \leq
\mu^{m-1}(1-\mu)$ so 
\begin{align*}
\gamma^{m-1} & \leq \mu^{m-1}\frac{1-\mu}{1-\gamma} \text{ and } \\
\gamma^m - \mu^m & \leq \frac{\gamma \mu^{m-1}(1-\mu)
 -  \mu^m(1-\gamma)}{1-\gamma} \\
& =   \frac{\mu^{m-1}(\gamma - \mu)}{1-\gamma} \\
& \leq   \frac{\gamma - \mu}{1-\gamma}, 
\end{align*}

\noindent
If instead of starting with a point $u \in \CD_{\mu} $ we had started
with a point $v \in \CD_{\gamma} $ we could have chosen  $u \in
\CD_{\mu} $ and made exactly the same computations as above.  This
symmetry in the choice of $u$ and $v$ 
shows we that the Hausdorff distance between  $\CD_{\mu} \text{ and
  } \CD_{\gamma}  $ is at most $2\frac{\gamma - \mu}{1 - \gamma},$ so
the Gromov--Hausdorff distance between $\CC_{\mu} \text{ and }
\CC_{\gamma} $ is at most this number too, and
the proposition follows.

\end{proof}

\subsection{ A compact metric space which contains compact Cantor sets of
  any Hausdorff dimension.}
The title for this subsection indicates it's content. By  the simple
definition $\CE_{\gamma}  := (1-\gamma)\CC_{\gamma}$ we define a
compact subset of $\ell^{1}(\IN, \IR) $ which is bi-Lipschitz
equivalent to $\CC_{\gamma} $. The closure of the union of all these
spaces, denoted $\CE$, is a compact space, which contains Cantor sets
of any dimension. In order to state the result a bit more precise we
define $e_1$ to be the unit vector in $\ell^{1}(\IN, \IR) $ which is
the first basis vector, i. e. the coordinates of $e_1$ are given by 
$e_1(n) = \delta_{1n}$. Further we define a subset $\CF $ of
$\ell^{1}(\IN, \IR) $ by $\CF = \{x \in \ell^{1}(\IN, \IR)\, : \, 0
\leq x(n) \leq 4/(n+1)^2 \, \}$. The set $\CF$ is a compact subset of
$\ell^{1}(\IN, \IR) $  since for elements  in $\CF$
we will have the following uniform estimates of the norm of the
tails of the elements
\begin{displaymath}
\forall x \in \CF\, \forall k \in \IN: \; \sum_{n=k}^{\infty} |x(n)|
\leq \frac{4}{k}.
\end{displaymath}

\begin{Theorem}
The space $\CE$ is a compact subset of $\CF $ of infinite Hausdorff
dimension. It contains closed Cantor sets of any positive Hausdorff
dimension and it can be described as the union
  \begin{displaymath}
\CE \, = \, {e_1} \cup \{ \underset{\gamma \in
  ]0,1[}{\cup}(1-\gamma)\CC_{\gamma} \} = \rm{closure}(\underset{\gamma \in
  ]0,1[}{\cup}\CE_{\gamma}).
\end{displaymath}
\end{Theorem}
 \begin{proof}
In order to prove that $\CE$ is contained in $ \CF$ we fix a $\gamma
$ in $]0,1[$ and look at an element  $x$ in $(1-\gamma)\CC_{\gamma}$.
For each $n$ in $\IN$  we have $ x(n)$ is in the two-point set  $ \{0,
(1- \gamma)^2\gamma^{n-1}\}$, so we must
 investigate the function $g_n(t) := (1-t)^2t^{n-1}$ on the unit
interval in order to find an upper estimate for it's maximum on this
interval.
 For $n = 1$ the maximal value of $g_1$ is $1 = 4/(1+1)^2$. 
For $n > 1$, elementary calculus yields that $g_n$ has maximal value
in the point $(n-1)/(n+1)$ and at this point
\begin{displaymath}
g_n(\frac{n-1}{n+1} )\, = \,( \frac{n-1}{n+1})^{n-1}(1- \frac{n-1}{n+1})^2
\, < \, (\frac{2}{n+1})^2,
\end{displaymath}
and $\CE$ is contained in $\CF$.

\noindent
In order to show that $e_1$ belongs to $\CE$ and that it is the only
point which has to be added to the union of the sets
$(1-\gamma)\CC_{\gamma}$, when closing up, we first remark that for
any $\gamma $ in $]0,1[$ $(1-\gamma)^2e_1 $ belongs to
$(1-\gamma)\CC_{\gamma}$, so $e_1$ must belong to $\CE$.
Suppose now that $x$ is a point in $\CE$.  Then there
exist  sequences $(\gamma_i)$ and $(x^i)$ such that for each $i$,
  $0 < \gamma_i < 1$
, $x^i $ is in $(1-\gamma_i)\CC_{\gamma_i}$ and $(x^i) $ converges to 
$x$. The sequence $(\gamma_i)$ is bounded so it has a convergent
subsequence and we may as well assume that the
sequence $(\gamma_i)$ is convergent with limit say $\xi$ in $[0,1]$.
If $\xi =1$ then $x =0$, but $0 $ is in all the sets $\CC_{\gamma}$,
so we may assume that $\xi < 1$. If $\xi = 0$ then $x = e_1$ and we
have dealt with this case already. We
can then assume that $0 < \xi < 1$. It is now quite elementary to
check that for  a coordinate, say $x(n)$   the
convergence of $(x^i)$  to $x$ implies that  
\begin{displaymath}
\forall n \in \IN: x(n) = 0 \text{ or }  x(n) =
(1-\xi)^2\xi^{n-1} 
\end{displaymath}
This means that $x$ is an element of $(1-\xi)\CC_{\xi}$ and the
theorem follows.
\end{proof}

\section{The uniform metric on the state space}

Let $\mathcal{A}$ be a unital C*-algebra. In this section we are studying the metric 
on the state space $\CS(\CA)$ which is induced by the norm on the dual space
$\CA^*$ of $\CA$.
For a  finite
dimensional C*-algebras the norm  topology and the weak*-topology
agree,
 so we looked for
a spectral triple for the algebra $\mathcal{M}_n$ of complex $ n \times n $
matrices such that the metric induced on the state space would be that
of the norm. As we mentioned in the introduction we considered 
a \textit{standard spectral triple}
given by $\mathcal{A}=\mathcal{M}_n$, $H=L^2(\mathcal{M}_n,
\frac{1}{n}\text{tr})$ and the operators in $\CA$ acting on $H$ by left
multiplication. The Dirac operator is then 
$D=T$ the selfadjoint unitary  operator on $H$ which
consists of transposing a matrix. 
It turned out that one can extend the idea behind the above spectral triple such that it is
possible, for any C*-algebra $\CA$, to construct a representation $\pi$ of  
$\CA$ on a Hilbert space $H$ such that there exists a projection $P
\in B(H)$ which has the property 
that the norm distance  on the state space is
recovered exactly if this projection $P$ plays the role of the Dirac
operator. This is the main result of this section. 
The proof is builded upon the following lemmas.
The first lemma is well known (\cite{Ri2}) and easy to prove.
\begin{Lemma} 
\label{unidist}
Let $A$ be a unital C*-algebra then for any two states $\phi, \psi$ on
$A$ 
\begin{displaymath} 
\|\; \phi - \psi \;\| \; = \; \sup \{ | \, (\phi - \psi) (a)\,| :\text{ }
 a = a^* \in \CA \; \rm{and} \; \underset{\alpha \in \IR}{\inf}\;\|\;a -
\alpha I\;\|\; \leq 1 \; \}.
\end{displaymath}
\end{Lemma}
\noindent
Further we state:
\begin{Lemma}
\label{unidirac}
Let $\CA$ be a unital C*-algebra and let $\rho$ denote a faithful
representation of $\CA$ on a Hilbert space $H$. Let $H_1$ denote the
Hilbert space tensor product $H_1 =H \otimes H$, $S$ the flip on $H_1$
given by $S (\xi \otimes \eta) = \eta \otimes \xi$ and $P$  the
projection $P = (I + S )/2$. Then the  representation $\pi$ of $\CA$ on 
$H_1$ given by the amplification $\pi(a) = \rho(a) \otimes I$ satisfies:
\begin{displaymath}
\forall a = a^* \in \CA: \quad \underset{\gamma \in \IR}{\inf}\;\|\;a - \gamma I\;\| \;= \; \|\;[P, \pi(a)]\; \| 
\end{displaymath}
\end{Lemma}
\begin{proof}
We will first transform the commutator slightly in order to ease the
computations, 
\begin{align*}
   \forall \gamma \in \IR \; \forall  a = a^* \in \CA :  \\
 \|[P, \pi(a)]\|
\;& = \;\frac{1}{2} \|[S, \pi(a)]\| \; = \; 
\frac{1}{2}\|S[S, \pi(a)]\| 
\\ &= \; \frac{1}{2} \|\pi(a) - S\pi(a)S\| 
  \\ &= \; \frac{1}{2}\|\rho(a) \otimes I - I \otimes 
\rho(a )\| 
\\ &= \; \frac{1}{2}\|\rho(a - \gamma I) \otimes I - I \otimes 
\rho(a - \gamma I)\| 
\end{align*} 

\noindent
From this series of identities it follows immediately that 
\begin{displaymath}
\forall a = a* \in \CA \, \forall \gamma \in \IR: \; \;  \|\;a -
\gamma I\;\|  \geq \; \|\;[P, \pi(a)]\; \|.
\end{displaymath}  

\noindent In order to show the  inequality in the opposite direction,
 for a certain $\gamma$, we use the series of identities again.
 Remark that by spectral theory it
follows that for $a= a^* \in \CA$ 
with spectrum contained in the smallest possible  interval
$[\alpha, \beta ] \subseteq \IR$ one has 
\begin{displaymath}
  \underset{\gamma \in \IR}{\inf}\;\|\;a - \gamma I\;\| \;= \; \|\;a
  - \frac{\alpha + \beta}{2}\; \| \; = \; \frac{\beta  - \alpha }{2}   
\end{displaymath}
Let $\eps > 0 $ and chose unit vectors $\xi , \eta \in H $ such
that $(\rho(a) \xi, \xi) \geq \beta - \eps $ and $(\rho(a) \eta,
\eta ) \leq \alpha  + \eps $. Then, $ \xi \otimes \eta $ is a
unit vector in $H_1$ and 
\begin{align*}
 \|\;[P, \pi(a)]\; \| & = 
\frac{1}{2} \|\rho(a) \otimes I - I \otimes 
\rho(a)\|\; \\ & \geq \;  \frac{1}{2}( (\rho(a) \otimes I - I \otimes \rho(a))
\xi \otimes \eta, \xi \otimes \eta )\; \\ & \geq \;  \frac{( \beta - \alpha )}{2} -
\eps \; \\ &= \;   \underset{\gamma \in \IR}{\inf}\;\|\;a -
\gamma I\;\|\; - \eps.
\end{align*}
The lemma follows.
\end{proof}

\noindent
We are  now ready to give the main result of this
section.
\begin{Theorem}
\label{unithm}
Let $\CA$ be a C*-algebra and $\rho$ a faithful non-degenerate
representation of $\CA$ on a Hilbert space $H$. Then there exists a
representation $\pi$ of $\CA$ on a Hilbert space $H_1$ which  
is an amplification of $\rho$, and a projection $P$ in
$B(H_1)$ such that for any pair of states $ \phi, \psi$ on $\CA$  
\begin{displaymath}
\|\; \phi - \psi \;\| \; = \; \sup \{ | \, (\phi - \psi) (a)\,| :\text{ }
 a = a^* \in \CA \; \rm{and} \; \|\; [P,\pi(a)]\; \| \; \leq 1 \; \}.
\end{displaymath}
If $H$ is separable and the commutant of $\rho(\CA)$ is a properly
infinite von Neumann algebra then $\pi = \rho$ is possible. 
If $\CA = M_n$ and $\rho$ is the  standard representation of $\CA$ on
$L^2( M_n, \frac{1}{n}\rm{tr})$ then $\pi = \rho$ is possible and
the projection $P = \frac{1}{2}( I + T )$ where $T$ is the
transposition on $M_n$ can be used.
\end{Theorem}
\begin{proof}
If $\CA$ has no unit then we add a unit in order to obtain a unital
C*-algebra $\tilde{\CA}$. It is well known that the state space of
$\CA$ embeds isometrically into  the state space of $\tilde{\CA}$. We
can then deduce the result for the non-unital case from the unital one
by remarking that both of the expressions 
\begin{displaymath}
| \, (\phi - \psi) (a)\,| \; \rm{and} \; \|\; [P,\pi(a)]\; \| 
\end{displaymath}
are left unchanged if $a$ is replaced by $(a - \alpha I )$.
\smallskip
\noindent Let us then assume that $\CA$ is unital. Then, by Lemma
\ref{unidirac} we can chose to amplify $\rho$ by the Hilbert dimension
of $H$, but less might do just as well. It all depends on the
multiplicity of the representation $\rho$, or rather whether the
commutant of $\rho(\CA)$ contains a subfactor isomorphic to $B(H)$.
In particular, this situation occurs if $H$ is separable and the
commutant is properly infinite.
\smallskip 
\noindent If $\CA = M_n$ and $\rho$ is the ``left regular representation''
of $\CA$ on $L^2( M_n, \rm{tr})$, then 
this Hilbert space is naturally identified with $\IC^n \otimes \IC^n$ via 
the mapping 
\begin{displaymath}
M_n(\IC) \ni a \to \sum_{i = 1}^n \sum_{j = 1}^n a_{ij} e_j
\otimes e_i,
\end{displaymath}
where the elements $e_i$ denote the elements of the
standard basis for $\IC^n$. From here it is easy to see that the flip
on the Hilbert space tensor product is nothing but the transposition
operator on $M_n$.
The  Lemma
\ref{unidirac} now applies directly for $\pi = \rho$ and the
projection $P = \frac{1}{2}( I + T)$, where $T$ is the transposition
operator on $M_n$. In the arguments above we have used the trace
rather than the trace-state as stated in the formulation of the
theorem, 
the reason
being that the identification of $M_n$ with $\IC^n \otimes \IC^n$
fits naturally with the trace.
\end{proof}


\begin{thebibliography}{99}

\bibitem[AC]{AC}
C. Antonescu, E. Christensen, {\em Metrics on group C*-algebras and a
  non-commutative Arzel\`{a}-Ascoli theorem}, J. Funct. Anal. 214 {\bf
  214} (2004), 247--259.


\bibitem[Br]{Br}
O. Bratteli, {\em Inductive limits of finite dimensional C*-algebras},
Trans. Amer. Math. Soc. {\bf 171} (1972), 195--234.



\bibitem[Ch]{Ch}
E. Christensen, {\em Subalgebras of a finite algebra}, Math. Anal, {\bf 243} (1979), 17--29.

\bibitem[Co1]{Co1}
A. Connes, {\em Classification of injective factors. Cases $II\sb{1},$
  $II\sb{\infty },$ $III\sb{\lambda },$ $\lambda \not=1$.},
  Ann. of Math. {\bf  104}  (1976),  73--115.  

\bibitem[Co2]{Co2}
A. Connes,  {\em Compact metric spaces, Fredholm modules, and
  hyperfiniteness},  Ergodic. Theor. and Dynam. Systems, {\bf 9} (1989), 207--220.

\bibitem[Co3]{Co3}
A. Connes,  {\em  Noncommutative Geometry}, Academic Press, San Diego,
1994.

\bibitem[Co4]{Co4}
A. Connes, {\em Unpublished notes on a Dirac operator associated to
  the Cantor subset of the unit interval.} 

\bibitem[Fa]{Fa}
K. J. Falconer, {\em Fractal Geometry, Mathematical Foundations and
  Applications}, Wiley 1990.  

\bibitem[FY]{FY}
H. Fukaishi, H. Yamaji, {\em An Elementary Construction of a Cantor
  set with Arbitrary Hausdorff Dimension }, RIMS Kokyuroku {\bf 1188},
 (2001), 86-95.

\bibitem[Gl]{Gl}
J. Glimm, {\em On a certain class of operator algebras}, Trans. Amer. Math. 
Soc., {\bf 95} (1960), 318--340.


\bibitem[Gr]{Gr}
M. Gromov, {\em Metric structures for Riemannian and non Riemannian
  spaces},
Birkhäuser, 1999. 

\bibitem[Haa]{Haa}
U. Haagerup,  {\em An example of a non nuclear C*-algebra which has the
  metric approximation property}, Invent. Math., {\bf
  50} (1979), 279--293.

 \bibitem[Jo]{Jo}
 P. Jolissaint, {\em Rapidly decreasing functions in reduced
  C*-algebras of groups}, Trans. Amer. Math. Soc., {\bf 317
  }(1990), 167--196.

\bibitem[KR]{KR}
R. V. Kadison, J. R. Ringrose,  {\em  Fundamentals of the theory of
  operator algebras}, Academic Press, New York, 1983.

\bibitem[Ke]{Ke}
J. Keesling, {\em Hausdorff dimension }, Topology Proc. {\bf
  11}(1986), 349 - 383. 

\bibitem[OR]{OR}
N. Ozawa, M. A. Rieffel, {\em Hyperbolic group C*-algebras and free 
product C*-algebras as compact quantum metric spaces},
Canad. J. Math., to appear,
arXiv:math. OA/0302310 v1, 2003.


\bibitem[Pav]{Pav}
B. Pavlovi\'{c},  {\em Defining metric spaces via operators from unital
  C*-algebras}, Pacific J. Math. {\bf 186} (1998), 285--313.


\bibitem[Ri1]{Ri1}
M. A. Rieffel, {\em Comments concerning non-commutative metrics}, 
Talk given at AMS Special Session, Texas A$\&$M, October 1993.

\bibitem[Ri2]{Ri2}
M.A. Rieffel  {\em Metrics on states from actions of compact groups},
Doc. Math., {\bf 3} (1998), 215--229.

\bibitem[Ri3]{Ri3}
M.A. Rieffel  {\em Metrics on state spaces},
Doc. Math. {\bf 4} (1999), 559--600.

\bibitem[Ri4]{Ri4}
M.A. Rieffel  {\em Group C*-algebras as compact quantum metric spaces},
Doc. Math. {\bf 7} (2002), 605--651.

\bibitem[Ri5]{Ri5}
M.A. Rieffel  {\em Compact Quantum Metric Spaces},
arXiv: math. OA/0308207 v1, 2003.

\bibitem [Ri6]{Ri6}
M.A. Rieffel  {\em Gromov--Hausdorff distance for quantum metric spaces,
}  Mem. Amer. Math. Soc.  {\bf 168 }  (2004),  no. 796, 1--65. 

\bibitem [Vo]{Vo}
D. V. Voiculescu  {\em On the existence of quasicentral approximate
  units relative to normed ideals. Part I},
J. Funct. Anal. {\bf 91} (1990), 1--36.


\end{thebibliography}
\end{document}